\documentclass[12pt,a4paper]{article}

\usepackage{amssymb,amsmath}
\usepackage[dvips]{graphicx}
\newtheorem{df}{Definition}[section]
\newtheorem{theo}[df]{Theorem}
\newtheorem{cor}[df]{Corollary}
\newtheorem{pr}[df]{Proposition}
\newtheorem{lem}[df]{Lemma}
\newtheorem{rem}[df]{Remark}
\newcommand{\qed}{\quad$\square$\par}
\newcommand{\dis}{\displaystyle}
\newcommand{\e}{\varepsilon}

\newcommand{\s}{\sigma}
\newcommand{\diag}{\mbox{diag}}

\newcommand{\cA}{{\cal A}}
\newcommand{\cB}{{\cal B}}

\newcommand{\cI}{{\cal I}}
\newcommand{\cL}{{\cal L}}

\newcommand{\cQ}{{\cal Q}}
\newcommand{\cT}{{\cal T}}

\newcommand{\bfi}{{\bf i}\hspace{1pt}}
\newcommand{\mR}{{\mathbb R}}

\newcommand{\br}{\mbox{\boldmath $\rho$}}
\newcommand{\bm}{\mbox{\boldmath $\mu$}}
\newcommand{\ba}{\mbox{\boldmath $a$}}
\newcommand{\bu}{\mbox{\boldmath $u$}}
\newcommand{\bv}{\mbox{\boldmath $v$}}

\newcommand{\bF}{\mbox{\boldmath $F$}}
\newcommand{\bG}{\mbox{\boldmath $G$}}
\newcommand{\bzero}{\mbox{\boldmath $0$}}
\newcommand{\bw}{\mbox{\boldmath $w$}}
\newcommand{\bV}{\mbox{\boldmath $V$}}
\renewcommand{\phi}{\varphi}
\newcommand{\bphi}{\mbox{\boldmath $\varphi$}}
\newcommand{\bkappa}{\mbox{\boldmath $\kappa$}}
\renewcommand{\tilde}{\widetilde}
\renewcommand{\bar}{\overline}
\makeatletter

\@addtoreset{equation}{section}
\makeatother
\begin{document}
\title{{\bf Nonlinear stability of stationary solutions
for curvature flow with triple junction}}
\author{Harald Garcke, Yoshihito Kohsaka, and Daniel \v{S}ev\v{c}ovi\v{c}}
\date{}
\maketitle

\pagestyle{myheadings}
\markright{STABILITY FOR CURVATURE FLOW WITH TRIPLE JUNCTION}

\noindent
{\bf Abstract.} 
In this paper we analyze the motion of a network of three planar curves 
with a speed proportional to the curvature of the arcs, having 
perpendicular intersections with the outer boundary and a common intersection 
at a triple junction.  As a main result we show that 
a linear stability criterion due to Ikota and Yanagida \cite{IY} is 
also sufficient for nonlinear stability. We also prove local and 
global existence of classical smooth solutions as well as various energy 
estimates. Finally, we prove exponential stabilization of an evolving 
network starting from the vicinity of a linearly stable stationary network.


\section{Introduction}

The motion of curves under the   curvature flow has been widely studied
in the past \cite{GH,Grayson,Ang}. Less is known about the evolution
of networks under the curvature flow \cite{BR,IY,Schnetal}. 
In this case the arcs in
the network evolve in the normal direction with a speed proportional to
the curvature of the arcs. At intersections with an outer boundary and
at triple junctions boundary conditions have to hold. At the outer
boundary one can prescribe the position (see \cite{KL,MNT}), or
the angle with the outer boundary \cite{BR, IY}. At the triple
junction Young's law, a force balance, leads to angle conditions. In
this paper we are interested in the stability of stationary solutions
to the  curvature flow with a triple junction when we prescribe
the natural angle condition of $90^\circ$ at the outer boundary. For
this case a linear stability criterion has been derived by Ikota and
Yanagida \cite{IY} (see also \cite{IY2}). We will demonstrate here
that this criterion also leads to nonlinear stability.

\begin{figure}
\centerline{
\includegraphics[width=5truecm]{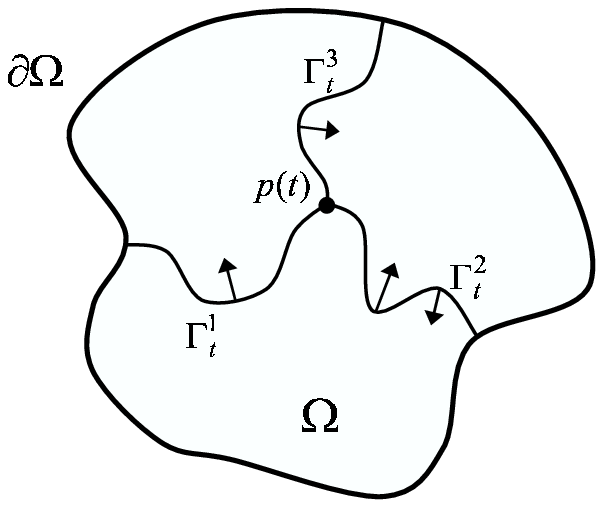}
\hglue 1truecm
\includegraphics[width=5truecm]{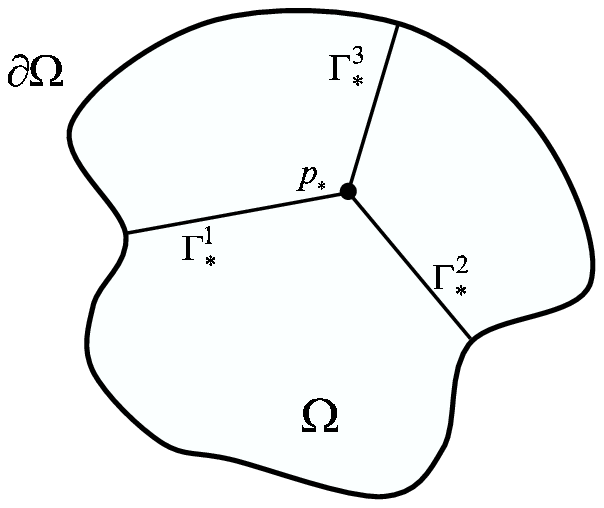}
}
\caption{
The curvature driven flow $\Gamma_t$ of a network with a triple junction 
at $p(t)$ (left) and a steady state $\Gamma_*$ (right).}
\label{fig1}
\end{figure}

We now specify the problem in detail. Let $\Omega$ be a bounded domain 
in $\mathbb{R}^2$ with $C^3$-boundary $\partial\Omega$. 
We introduce a $C^3$-function $\psi:\mathbb{R}^2 \to\mathbb{R}$ 
with $\nabla\psi(x)\ne 0$ if $\psi(x)=0$ such that
$$
\Omega=\{x\in\mathbb{R}^2\mid\psi(x)<0\},\qquad
\partial\Omega=\{x\in\mathbb{R}^2\mid \psi(x)=0\}.
$$
We search families of curves $\Gamma^1_t$, $\Gamma^2_t$, and $\Gamma^3_t$ 
which are parameterized by time $t$ and which are contained in $\Omega$. The
three curves are supposed to meet at a triple junction $p(t)\in\Omega$
at their one end point and at the other end point they are required
to intersect with $\partial \Omega$, see Figure~\ref{fig1}.
We require for $i=1,2,3$
\begin{align}
\hspace*{3cm}
&\beta^iV^i=\gamma^i\kappa^i \hspace*{-2.2cm}
&\mbox{on} \hspace*{0.4cm}
&\Gamma^i_t\,,& \label{flow} \\
&\dis\sum^3_{i=1}\gamma^iT^i=0 \hspace*{-2.2cm}
&\mbox{at} \hspace*{0.4cm}
&p(t),& \label{YL} \\
&\Gamma^i_t\bot\partial\Omega \hspace*{-2.2cm}
&\mbox{at} \hspace*{0.4cm}
&\Gamma^i_t\cap\partial\Omega.& \label{outbc}
\end{align}
Here $V^i$ and $\kappa^i$ are the normal velocity and curvature  of $\Gamma^i_t,$ respectively. The constants $\beta^i$ and $\gamma^i$ 
are given physical parameters and $T^i$ are unit tangents to the curve
which are chosen such that they point away from the triple junction. 

Equation (\ref{YL}) is a force balance and one can solve for
the $T^i$'s if the condition
\begin{equation*}
\gamma^i+\gamma^j\leq \gamma^k 
\quad\mbox{for all}\,\ \{i,j,k\}\,\ \mbox{mutually different}\,,
\end{equation*}
is fulfilled. In the following we assume strict inequalities and an
argument as in Bronsard and Reitich \cite{BR} gives that the angles 
$\theta^i$ between the tangents $T_j$ and $T_k$ fulfill 
\begin{equation*}
\frac{\sin\theta^1}{\gamma^1}=\frac{\sin\theta^2}{\gamma^2}
=\frac{\sin\theta^3}{\gamma^3}
\end{equation*}
with $0<\theta^i<\pi$ ($i=1,2,3$) and $\theta^1+\theta^2+\theta^3=2\pi$. 
Existence of solutions to the evolution problem (\ref{flow})-(\ref{outbc})
has been shown by Bronsard and Reitich \cite{BR}. We will show later that 
the energy functional 
\begin{equation*}
E[\Gamma_t] = \sum^3_{i=1} \gamma^i L[\Gamma^i_t]
\end{equation*}
where $\Gamma_t=\bigcup^3_{i=1}\Gamma^i_t$ and $L[\Gamma^i_t]$ is 
the length of $\Gamma^i_t$, is a Ljapunov functional. 
The constants $\gamma^i$ can be interpreted as surface free energy 
densities (surface tensions) and the functional $E$ is the total 
free energy of the systems. Sternberg and Zeimer \cite{SZ} showed 
the existence of isolated local minimizers to $E$, which can be 
interpreted as solutions to a partitioning problem of two dimensional 
domains into three subdomains having (locally) least interfacial area.

The paper is organized as follows. In the next section we present 
a way how to parameterize the problem. We derive a nonlinear nonlocal 
system of parabolic equations governing the evolution of curves driven 
by  curvature. By means of the semi-group theory due to Lunardi \cite{Lu} 
we prove local existence of a classical solution. 
Section~\ref{lin} is devoted to the rigorous derivation of 
the linearized system of equations. We recall the result of 
Yanagida and Ikota stating an explicit condition for linearized stability 
of the governing system of equations. In Section~\ref{uniqueness} 
we provide a usefull result guaranteeing local uniqueness of 
a stationary solution proved by the inverse function theorem 
and the result is to our knowledge the first result in this direction for
networks. As a byproduct we also obtain an important bound for
the displacement of the network in terms of the curvature.  We proceed by deriving 
useful geometric equations for the curvature and other geometric 
quantities in  Section \ref{goveq}. Using the linearized stability 
criterion we show how  to derive a priori estimates for
Sobolev norms of the solution.
These geometric equations are then used in order to prove usefull 
bounds for a solution. With the help of these energy type estimates 
we prove global existence of a classical solution. In the final 
section~\ref{stab} we prove exponential 
stability of the stationary solution.

\section{Parameterization and local existence}\label{sec:param}

\begin{figure}
\centerline{
\includegraphics[angle=270]{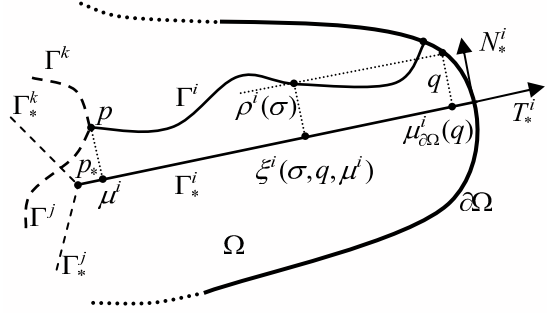}
}
\caption{
Description of the local parametrization of the curve $\Gamma^i$.}
\label{fig2}
\end{figure}

We consider line segments $\Gamma^1_\ast$, $\Gamma^2_\ast$ and 
$\Gamma^3_\ast$ meeting the outer boundary with an angle of $90^\circ$ 
at their one end point and having without loss of generality 
$p_\ast=(0,0)^T$ as their common other end point where we assume 
that (\ref{YL}) holds. 
Then we define an arc-length parameterization of $\Gamma^i_\ast$ 
($i=1,2,3$) as
\begin{equation*}
\Gamma^i_\ast = \{\Phi^i_\ast(\sigma) \mid \sigma\in [0,l^i]\}
\end{equation*}
with $\Phi^i_\ast(0)=p_\ast = (0,0)^T$, $\Phi^i_\ast(l^i)\in\partial\Omega$. 
In particular, we obtain that $l^i$ is the length of $\Gamma^i_\ast$. 
Then we will extend $\Phi^i_\ast$ as an arc-length parameterization of 
the full line which contains $\Gamma^i_\ast$. 
We will now introduce a certain stretched coordinate system 
in order to allow for parameterizations of curves close to $\Gamma^i_\ast$ 
($i=1,2,3$) over fixed intervals $[0,l^i]$. 

Let $T^i_\ast$ be the unit tangent to $\Gamma^i_\ast$ pointing from 
the triple junction $p_\ast$ to the outer boundary and 
let $N^i_\ast=RT^i_\ast$ be a unit normal where $R$ is the anticlockwise rotation by $\pi/2$. We then define
\begin{equation*}
\mu^i=(p,T^i_\ast)_{\mR^2}, \quad 
\mu_{\partial\Omega}^i(q)=\max\{\,\sigma \mid 
\Phi^i_\ast(\sigma)+qN^i_\ast\in\overline{\Omega}\,\}\,.
\end{equation*}
We remark that the parameter $\mu^i$ allows for a tangential movement 
of the triple junction along $\Gamma_*^i$. We now set
\begin{equation*}
\Psi^i(\sigma,q,\mu^i)=\Phi_\ast^i(\xi^i(\sigma,q,\mu^i))+qN^i_\ast\,,
\end{equation*}
where
\begin{equation*}
\xi^i(\sigma,q,\mu^i)
=\mu^i+\frac{\sigma}{l^i}(\mu_{\partial\Omega}^i(q)-\mu^i)\,.
\end{equation*}
Note that $\xi^i(\sigma,0,0)=\sigma$ and $\xi^i(0,q,\mu^i)=\mu^i$. 

We now define the parameterization of curves 
$\Gamma=(\Gamma^1,\Gamma^2,\Gamma^3)$ close to 
$\Gamma_\ast=(\Gamma^1_\ast,\Gamma^2_\ast,\Gamma^3_\ast)$ 
having their triple junction at the point $p$ with the help of functions
\begin{equation*}
\rho^i:[0,l^i]\to\mathbb{R}
\end{equation*}
which fulfill the conditions
\begin{equation}
\rho^i(0)=(p,N_\ast^i)_{\mR^2} \quad (i=1,2,3)
\label{rho_at0}
\end{equation}
(see Fig.~\ref{fig2}). Set
\begin{equation}\label{param}
\Phi^i(\sigma)
=\Psi^i(\sigma,\rho^i(\sigma),\mu^i)\,,\,\,\sigma\in[0,l^i]\,.
\end{equation}
Then the functions $\Phi^i$ parameterize the curves $\Gamma^i$ 
in the neighborhood of $\Gamma_\ast$ as 
$\Gamma^i=\{\Phi^i(\sigma)\,|\,\sigma\in[0,l^i] \}$.
Since $\Phi^i_\ast(\mu^i)=\mu^iT_\ast^i$, we have 
$\Phi^i(0)=\mu^iT_\ast^i+\rho^i(0)N_\ast^i=p$, which implies that 
\begin{equation}\label{stick}
\mu^1T^1_\ast+\rho^1(0)N^1_\ast
=\mu^2T^2_\ast+\rho^2(0)N^2_\ast
=\mu^3T^3_\ast+\rho^3(0)N^3_\ast\,.
\end{equation}
By virtue of the definition of $\mu^i$, equation (\ref{rho_at0}), and Young's law 
\begin{equation*}
\sum^3_{i=1}\gamma^iT^i_\ast = 0, \quad
\sum^3_{i=1}\gamma^iN^i_\ast=0,
\end{equation*}
we are led to
\begin{equation*}
\sum^3_{i=1}\gamma^i\mu^i=0, \quad
\sum^3_{i=1}\gamma^i\rho^i(0)=0\,.
\end{equation*}
Furthermore, identities (\ref{stick}) and the angle conditions give the following lemma. 

\begin{lem} \label{lem:junction}
Let us define the matrix
\[
Q=-\frac1{1-c^1c^2c^3}
\left(\begin{array}{ccc}
c^3c^1s^2&s^3&c^3s^1 \\
c^1s^2&c^1c^2s^3&s^1 \\
s^2&c^2s^3&c^2c^3s^1
\end{array}\right),
\]
where $c^i=\cos\theta^i$, $s^i=\sin\theta^i$. Then, 
for $\br=(\rho^1,\rho^2,\rho^3)$, and $\bm=(\mu^1,\mu^2,\mu^3)$, 
it holds $\bm^T=Q\br^T(0)$.
\end{lem}
{\it Proof.}
It follows from (\ref{stick}) that 
$$
\mu^i=\mu^j(T_*^j,T_*^i)_{\mR^2}+\rho^j(0)(N_*^j,T_*^i)_{\mR^2}
$$
for $(i,j)=(1,2),\,(2,3),\,(3,1)$. By the angle condition, we have 
$$
(T_*^i,T_*^j)_{\mR^2}=\cos\theta^k, \quad 
(T_*^i,N_*^j)_{\mR^2}=\cos(\theta^k+\pi/2)=-\sin\theta^k
$$
for $(i,j,k)=(1,2,3),\,(2,3,1),\,(3,1,2)$. This implies 
$\mu^i-c^k\mu^j=-s^k\rho^j(0)$, so that we are led to 
$$
\left(\begin{array}{ccc}
-c^2&0&1 \\ 1&-c^3&0 \\ 0&1&-c^1
\end{array}\right)
\left(\begin{array}{c}
\mu^1 \\ \mu^2 \\ \mu^3
\end{array}\right)
=\left(\begin{array}{ccc}
-s^2&0&0 \\ 0&-s^3&0 \\ 0&0&-s^1
\end{array}\right)
\left(\begin{array}{c}
\rho^1(0) \\ \rho^2(0) \\ \rho^3(0)
\end{array}\right).
$$
Then we obtain
\[
\left(\begin{array}{c}
\mu^1 \\ \mu^2 \\ \mu^3
\end{array}\right)
=\frac{-1}{1-c^1c^2c^3}
\left(\begin{array}{ccc}
c^3c^1&1&c^3 \\ c^1&c^1c^2&1 \\ 1&c^2&c^2c^3
\end{array}\right)
\left(\begin{array}{ccc}
s^2&0&0 \\ 0&s^3&0 \\ 0&0&s^1
\end{array}\right)
\left(\begin{array}{c}
\rho^1(0) \\ \rho^2(0) \\ \rho^3(0)
\end{array}\right) 
= Q \left(\begin{array}{c}
\rho^1(0) \\ \rho^2(0) \\ \rho^3(0)
\end{array}\right) 
\]
which completes the proof. \qed

\smallskip

We now consider evolving curves
\begin{equation*}
\Gamma^i(t)=\{\Phi^i(\sigma,t)\mid\sigma\in[0,l^i]\},
\end{equation*}
where $\Phi^i(\cdot,t)$ are defined as in (\ref{param}) such that 
$\rho^i(\cdot,t)$ ($i=1,2,3$) satisfy (\ref{rho_at0}). 
We formulate the curvature flow for a network with the help of
these parameterizations. For that purpose, the following quantities 
are needed
\begin{eqnarray}
&&T^i=\frac{1}{|\Phi^i_\sigma|}\Phi^i_\sigma, \quad
N^i=\frac{1}{|\Phi^i_\sigma|}R\Phi^i_\sigma\,,\nonumber \\[0.1cm]
&&V^i=(\Phi^i_t,N^i)_{\mR^2}\,,\ \ 
\kappa^i=\bigl(\,\frac{1}{|\Phi^i_\sigma|}T^i_{\s},
N^i\bigr)_{\mR^2}\,, \ \ i=1,2,3\,. \label{kappai}
\end{eqnarray}
Then we obtain the following formulation for the curvature flow of 
a network:
\begin{eqnarray}
&&\beta^i(\Phi^i_t,N^i)_{\mR^2}
=\gamma^i\bigl(\,\frac{1}{|\Phi^i_\sigma|}T^i_{\sigma},
N^i\bigr)_{\mR^2}\,,  \ \ i=1,2,3\,,\label{p1} \\
&&\sum^3_{i=1}\gamma^i\rho^i=0 \quad\mbox{at}\quad 
\sigma=0\,,\label{p2} \\[0.1cm]
&&(T^1,T^2)_{\mR^2}=\cos\theta^3,\quad   
(T^3,T^1)_{\mR^2}=\cos\theta^2 \quad\mbox{at}\quad
\sigma=0\,,\label{p3_2} \\[0.2cm]
&&(N^i,\nabla\psi(\Phi^i))_{\mR^2}=0 \quad\mbox{at}\quad
\sigma=l^i\,,  \ \ i=1,2,3\,.\label{p4}
\end{eqnarray}
Note that the conditions $\psi(\Phi^i(l^i,t))=0$ ($i=1,2,3$) and
$\Phi^1(0,t)=\Phi^2(0,t)=\Phi^3(0,t)$ are always fulfilled 
with our choice of the parameterizations. 
We can formulate the problem in terms of $(\rho^1,\rho^2,\rho^3)$ 
and obtain a system of three second order parabolic equations 
where each equation is defined on a different spatial interval $[0,l^i]$. 
We obtain one boundary condition at $\sigma=l^i$ for the $i$-th equation 
and the three equations are coupled through the three boundary conditions 
at the triple junction. 

Let us derive the form of the nonlinear system for $\rho^i(\s,t)$ 
($i=1,2,3$). Set $\cI^i=[0,l^i]$ and $\cQ^i_{t_0,t_1}
=\cI^i\times(t_0,t_1]$. Equations (\ref{p1}) give
\begin{equation}
\rho^i_t=L^i(\rho^i,\rho^i_\s,\mu^i)
\kappa^i(\rho^i,\rho^i_\s,\rho^i_{\s\s},\mu^i)
+\Lambda^i(\rho^i,\rho^i_\s,\mu^i)\mu^i_t
\quad\mbox{for}\,\ (\s,t)\in\cQ^i_{0,T},
\label{rho_eq}
\end{equation}
where $L^i(\rho^i,\rho^i_\s,\mu^i)$ and 
$\Lambda^i(\rho^i,\rho^i_\s,\mu^i)$ are
\begin{eqnarray*}
&&L^i(\rho^i,\rho^i_\s,\mu^i)
=\frac{\gamma^i}{\beta^i(\Psi^i_q,R\Psi^i_{\s})_{\mR^2}}\hspace*{1pt}
J^i(\rho^i,\rho^i_\s,\mu^i), \\[0.1cm]
&&\Lambda^i(\rho^i,\rho^i_\s,\mu^i)
=-\frac1{(\Psi^i_q,R\Psi^i_{\s})_{\mR^2}}
\bigl\{(\Psi^i_{\mu},R\Psi^i_{\s})_{\mR^2}
+(\Psi^i_{\mu},R\Psi^i_q)_{\mR^2}\rho^i_{\s}\bigr\}
\end{eqnarray*}
with the notation $J^i=|\Phi^i_\s|$, and the curvature 
$\kappa^i=\kappa^i(\rho^i,\rho^i_\s,\rho^i_{\s\s},\mu^i)$ 
is represented as
\begin{eqnarray*}
&&\kappa^i(\rho^i,\rho^i_\s,\rho^i_{\s\s},\mu^i) \\[0.1cm]
&&=\frac1{\{J^i(\rho^i,\rho^i_\s,\mu^i)\}^3}\Bigl[
(\Psi^i_q,R\Psi^i_{\s})_{{\mathbb R}^2}\rho^i_{\s\s}
+\bigl\{2(\Psi^i_{\s q},R\Psi^i_{\s})_{{\mathbb R}^2}
+(\Psi^i_{\s\s},R\Psi^i_q)_{{\mathbb R}^2}\bigr\}\rho^i_{\s} 
\nonumber \\[0.3cm]
&&\hspace*{0.5cm}
+\bigl\{(\Psi^i_{qq},R\Psi^i_{\s})_{{\mathbb R}^2}
+2(\Psi^i_{\s q},R\Psi^i_q)_{{\mathbb R}^2}
+(\Psi^i_{qq},R\Psi^i_q)_{{\mathbb R}^2}\rho^i_{\s}\bigr\}(\rho^i_{\s})^2 
+(\Psi^i_{\s\s},R\Psi^i_{\s})_{{\mathbb R}^2}\Bigr].
\end{eqnarray*}
By virtue of Lemma \ref{lem:junction} and (\ref{rho_eq}), we have
$$
\left(\begin{array}{c}
\mu^1_t \\ \mu^2_t \\ \mu^3_t
\end{array}\right)
=Q\hspace{1pt}\cT^0\hspace{-1pt}\left(\begin{array}{c}
\rho^1_t \\ \rho^2_t \\ \rho^3_t
\end{array}\right)
=Q\Bigl[\cT^0\hspace{-1pt}M(\br,\br_{\s},\bm)\Bigr]^{-1}
\cT^0\hspace{-1pt}\left(\begin{array}{c}
L^1(\bu^1)\kappa^1(\bu^1) \\ 
L^2(\bu^2)\kappa^2(\bu^2) \\ 
L^3(\bu^3)\kappa^3(\bu^3)
\end{array}\right), 
$$
where $Q$ is the matrix as in Lemma \ref{lem:junction}, 
$\cT^0$ is the trace operator to $\s=0$, i.e. 
$\cT^0\hspace{-1pt}f=f\big|_{\s=0}$, and 
$M(\br,\br_{\s},\bm)$ is the matrix 
\[ 
M(\br,\br_{\s},\bm) 
=\mbox{Id} - \diag(\Lambda^1(\bu^1), \Lambda^2(\bu^2), \Lambda^3(\bu^3)) Q
\]
with the notation $\bu^i=(\rho^i,\rho^i_\s,\mu^i)$. 

\begin{rem}
The matrix $M(\br,\br_{\s},\bm)$ is invertible 
provided that 
\begin{equation}
\sum_{i=1}^3\Bigl(\,\sup_{t\in[0,T]}\|\rho^i(\cdot,t)\|_{C^1(\cI^i)}
+\sup_{t\in[0,T]}|\mu^i(t)|\,\Bigr)<\delta_0
\label{c-inverse}
\end{equation}
for some $\delta_0>0$. Indeed, we have 
$\det M(\br,\br_{\s},\bm)
=d\bigl\{-1+(c^2-\Lambda^1(\bu^1)s^2)(c^3-\Lambda^2(\bu^2)s^3)
(c^1-\Lambda^3(\bu^3)s^1)\bigr\}$, where $d=-1/(1-c^1c^2c^3)$. 
Then, $\Lambda^i(\bzero)=0$ ($i=1,2,3$) imply that   
$\det M(\bzero,\bzero,\bzero)=d(-1+c^1c^2c^3)=1\ne0$.
Since $\det M(\br,\br_{\s},\bm)$ is continuous with respect to 
$\br$, $\br_{\s}$, and $\bm$, we can conclude that for $\e<1$ 
there exists a $\delta_0>0$ such that 
$\det M(\br,\br_{\s},\bm)>1-\e>0$ provided that (\ref{c-inverse}) holds. 
\end{rem}

\noindent
As a consequence, we are led to the following nonlinear nonlocal partial 
differential equations for $\rho^i(\s,t)$ ($i=1,2,3$):
\begin{eqnarray*}
\rho^i_t
&=&a^i(\rho^i,\rho^i_\s,\mu^i)\rho^i_{\s\s}
+\Lambda^i(\rho^i,\rho^i_\s,\mu^i)
\sum_{j=1}^3a_1^{ij}(\cT^0\hspace{-1pt}\br,\cT^0\hspace{-1pt}\br_\s,\bm)
\cT^0\hspace{-1pt}\rho^j_{\s\s} \\
&&+f^i(\rho^i,\partial_\s\rho^i,\cT^0\hspace{-1pt}\br,
\cT^0\hspace{-1pt}\br_\s,\bm) 
\quad\mbox{for}\,\ (\s,t)\in\cQ^i_{0,T}
\end{eqnarray*}
where $a^i(\rho^i,\rho^i_\s,\mu^i)
=\gamma^i/\beta^i\{J^i(\rho^i,\rho^i_\s,\mu^i)\}^2$ 
and $a_1^{ij}(\cT^0\hspace{-1pt}\br,\cT^0\hspace{-1pt}\br_\s,\bm)$ 
is the $(i,j)$-component of the matrix 
\[
\ba_1(\cT^0\hspace{-1pt}\br,\cT^0\hspace{-1pt}\br_\s,\bm)
=Q\Bigl[\cT^0\hspace{-1pt}M(\br,\br_{\s},\bm)\Bigr]^{-1}
\diag(\cT^0\hspace{-1pt}a^1(\bu^1), \cT^0\hspace{-1pt}a^2(\bu^2), 
\cT^0\hspace{-1pt}a^3(\bu^3))\,.
\]
Furthermore, $f^i$ is a smooth function in $\mR$ which is evaluated at 
lower order terms. 
Then, recalling the boundary conditions (\ref{p2})-(\ref{p4}) and 
Lemma \ref{lem:junction}, we have the following nonlinear system:
\begin{equation}
\left\{\begin{array}{l}
\rho^i_t
=a^i(\rho^i,\rho^i_\s,\mu^i)\rho^i_{\s\s}
+\Lambda^i(\rho^i,\rho^i_\s,\mu^i)
\dis\sum_{j=1}^3a_1^{ij}(\cT^0\hspace{-1pt}\br,
\cT^0\hspace{-1pt}\br_\s,\bm)
\cT^0\hspace{-1pt}\rho^j_{\s\s} \\[0.6cm]
\hspace*{0.85cm}
+f^i(\rho^i,\rho^i_\s,\cT^0\hspace{-1pt}\br,\cT^0\hspace{-1pt}\br_\s,\bm)
\quad\mbox{for}\,\ (\s,t)\in \cQ^i_{0,T}, \\[0.2cm]
\dis\sum_{i=1}^3\gamma^i \rho^i=0, \quad  
g^{12}(\bu^{12})=0, \quad g^{13}(\bu^{13})=0
\quad\mbox{at}\,\ \s=0, \\[0.6cm]
b_{\partial\Omega}^i(\rho^i,\mu^i)\rho^i_{\s}
+g_{\partial\Omega}^i(\rho^i,\mu^i)=0
\quad\mbox{at}\,\ \s=l^i \quad (i=1,2,3), \\[0.3cm]
\bm^T=Q\hspace{1pt}\cT^0\hspace{-1pt}\br^T
\quad\mbox{for}\,\ t\in(0,T]
\end{array}\right.
\label{NL}
\end{equation}
where $\bu^{1j}=(\rho^1,\rho^j,\rho^1_\s,\rho^j_\s,\mu^1,\mu^j)$ 
($j=2,3$) and 
\begin{eqnarray*}
&&g^{12}(\bu^{12})
=(\Psi^1_{\s},\Psi^2_{\s})_{\mR^2}
+(\Psi^1_{\s},\Psi^2_q)_{\mR^2}\rho^2_{\s}
+(\Psi^1_q,\Psi^2_{\s})_{\mR^2}\rho^1_{\s}
+(\Psi^1_q,\Psi^2_q)_{\mR^2}\rho^1_{\s}\rho^2_{\s} \\[0.15cm]
&&\hspace*{2.1cm}
-J^1(\rho^1,\rho^1_\s,\mu^1)J^2(\rho^2,\rho^2_\s,\mu^2)\cos\theta^3, 
\\[0.15cm]
&&g^{13}(\bu^{13})
=(\Psi^3_{\s},\Psi^1_{\s})_{\mR^2}
+(\Psi^3_{\s},\Psi^1_q)_{\mR^2}\rho^1_{\s}
+(\Psi^3_q,\Psi^1_{\s})_{\mR^2}\rho^3_{\s}
+(\Psi^3_q,\Psi^1_q)_{\mR^2}\rho^3_{\s}\rho^1_{\s} \\[0.15cm]
&&\hspace*{2.1cm}
-J^3(\rho^3,\rho^3_\s,\mu^3)J^1(\rho^1,\rho^1_\s,\mu^1)\cos\theta^2, 
\\[0.15cm]
&&b_{\partial\Omega}^i(\rho^i,\mu^i)
=(R\Psi^i_q,\nabla\psi(\Psi^i))_{\mR^2}, \quad
g_{\partial\Omega}^i(\rho^i,\mu^i)
=-(R\Psi^i_{\s},\nabla\psi(\Psi^i))_{\mR^2}.
\end{eqnarray*}
Now we are ready to state a local existence result. 

\begin{theo}[Local existence]\label{th:l-exist}
Let $\alpha\in(0,1)$ and let us assume that 
$\rho_0^i\in C^{2+\alpha}(\cI^i)$ and $\mu_0^i$ ($i=1,2,3$) with 
sufficiently small norms $\|\rho_0^i\|_{C^{1+\alpha}(\cI^i)}$ and $|\mu_0^i|$ 
fulfill the compatibility conditions
$$
\left\{\begin{array}{l}
\dis\sum_{i=1}^3\gamma^i\rho^i_0=0, \quad 
g^{12}(\bu^{12}_0)=0, \quad g^{13}(\bu^{13}_0)=0
\quad\mbox{at}\,\ \s=0, \\[0.4cm]
b_{\partial\Omega}^i(\rho^i_0,\mu^i_0)\rho^i_{0,\s}
+g_{\partial\Omega}^i(\rho^i_0,\mu^i_0)=0
\quad\mbox{at}\,\ \s=l^i \quad (i=1,2,3),
\end{array}\right.
$$
where $\bu^{1j}_0=(\rho^1_0,\rho^j_0,\rho^1_{0,\s},\rho^j_{0,\s},
\mu^1_0,\mu^j_0)$ ($j=2,3$). 
Then there exists a 
$$
T_0=T_0\Bigl(1/\sum_{i=1}^3\|\rho_0^i\|_{C^{2+\alpha}(\cI^i)}\Bigr)>0\,,
$$
$T_0$ being an increasing function of its argument and 
such that the problem (\ref{NL}) with 
$(\rho^i(\cdot,0),\mu^i(0))=(\rho_0^i,\mu_0^i)$ ($i=1,2,3$) has 
a unique solution 
$$
(\rho^1,\rho^2,\rho^3,\mu^1,\mu^2,\mu^3)\in 
C^{2+\alpha,1}(\overline{\cQ^1_{0,T_0}})
\times C^{2+\alpha,1}(\overline{\cQ^2_{0,T_0}})
\times C^{2+\alpha,1}(\overline{\cQ^3_{0,T_0}})
\times \bigl[C^1[0,T_0]\bigr]^3
$$
satisfying (\ref{c-inverse}). 
\end{theo}

In order to prove Theorem \ref{th:l-exist} by using a contraction principle, 
we need some preparations which consist of three steps: 1) 
the linearization of (\ref{NL}) around the initial data; 2)
the verification of the complementary conditions for the linearized system; 
3) the derivation of suitable a priori estimate for solutions of 
the linearized system. 

\smallskip

{\it Step 1.} 
Let us derive the linearization of (\ref{NL}) around the initial data 
$\rho_0^i\in C^{2+\alpha}(\cI^i)$ and $\mu_0^i$ ($i=1,2,3$). 
First we define differential operators as
\begin{eqnarray*}
&&\cA_0=\diag(a^1(\bu_0^1), a^2(\bu_0^2), a^3(\bu_0^3)) 
\partial_\s^2, \\[0.1cm]
&&\cA_1=\diag(\Lambda^1(\bu_0^1),\Lambda^2(\bu_0^2),\Lambda^3(\bu_0^3)) 
\ba_1(\cT^0\hspace{-1pt}\br_0,\cT^0\hspace{-1pt}\br_{0,\s},\bm_0)
\cT^0\partial_\s^2,
\end{eqnarray*}
and also define, for given functions $(\bar{\rho}^i,\bar{\mu}^i)\in 
C^{2+\alpha,1}(\overline{\cQ^i_{0,T}})\times C^1[0,T]$ ($i=1,2,3$),
\begin{eqnarray*}
F^i(\s,t)
&=&\bigl\{a^i(\bar{\bu}^i)-a^i(\bu_0^i)\bigr\}\bar{\rho}^i_{\s\s} \\
&&+\sum_{j=1}^3\bigl\{
\Lambda^i(\bar{\bu}^i)
a_1^{ij}(\cT^0\hspace{-1pt}\bar{\br},\cT^0\hspace{-1pt}\bar{\br}_\s,
\bar{\bm})
-\Lambda^i(\bu_0^i)
a_1^{ij}(\cT^0\hspace{-1pt}\br_0,\cT^0\hspace{-1pt}\br_{0,\s},\bm_0)
\bigr\}\cT^0\hspace{-1pt}\bar{\rho}^j_{\s\s} \\[0.1cm]
&&+f^i(\bar{\rho}^i,\bar{\rho}^i_\s,
\cT^0\hspace{-1pt}\bar{\br},\cT^0\hspace{-1pt}\bar{\br}_\s,\bar{\bm}),
\end{eqnarray*}
where $\bu^i_0=(\rho_0^i,\rho^i_{0,\s},\mu_0^i)$, 
$\br_0=(\rho^1_0,\rho^2_0,\rho^3_0)$, 
$\bm_0=(\mu^1_0,\mu^2_0,\mu^3_0)$, 
$\bar{\bu}^i=(\bar{\rho}^i,\bar{\rho}^i_{\s},\bar{\mu}^i)$, 
$\bar{\br}=(\bar{\rho}^1,\bar{\rho}^2,\bar{\rho}^3)$, and 
$\bar{\bm}=(\bar{\mu}^1,\bar{\mu}^2,\bar{\mu}^3)$. 
Then, setting $\bF=(F^1,F^2,F^3)$, we have the linearization of 
the differential equation given as
$$
\partial_t\br^T=\cA_0\br^T+\cA_1\br^T+\bF^T(\s,t).
$$
Let us derive the linearization of the boundary conditions. 
For $\rho_0^i\in C^{2+\alpha}(\cI^i)$ and $\mu_0^i$ ($i=1,2,3$), 
we define differential operators as
$$
(B^{ki}(\cT^0\hspace{-1pt}\br_0,
\cT^0\hspace{-1pt}\br_{0,\s},\bm_0)\partial_\s)_{i=1,2,3} 
=
\left\{
\begin{matrix}
(\gamma^1,\gamma^2,\gamma^3)\hfill 
\,\mbox{for}&k=1, \\[0.2cm]
(b^{21}(\cT^0\bu^{12}_0)\partial_\s,b^{22}(\cT^0\bu^{12}_0)\partial_\s,0)
\,\mbox{for}&k=2, \\[0.2cm]
(b^{31}(\cT^0\bu^{13}_0)\partial_\s,0,b^{33}(\cT^0\bu^{13}_0)\partial_\s)
\,\mbox{for}&k=3.
\end{matrix}
\right.
$$

Here the components are  represented as follows:
\begin{eqnarray*}
b^{21}(\bu^{12})=(\Psi^1_q,\Psi^2_{\s})_{\mR^2}
+(\Psi^1_q,\Psi^2_q)_{\mR^2}\,\rho^2_{\s} 
-\bigl\{(\Psi^1_{\s},\Psi^1_q)_{\mR^2}
+|\Psi^1_q|^2\,\rho^1_{\s}\bigr\}
\frac{J^2(\bu^2)}{J^1(\bu^1)}\cos\theta^3, \\[0.2cm]
b^{22}(\bu^{12})=(\Psi^1_{\s},\Psi^2_q)_{\mR^2}
+(\Psi^1_q,\Psi^2_q)_{\mR^2}\,\rho^1_{\s} 
-\bigl\{(\Psi^2_{\s},\Psi^2_q)_{\mR^2}
+|\Psi^2_q|^2\,\rho^2_{\s}\bigr\}
\frac{J^1(\bu^1)}{J^2(\bu^2)}\cos\theta^3, \\[0.2cm]
b^{31}(\bu^{13})=(\Psi^3_{\s},\Psi^1_q)_{\mR^2}
+(\Psi^3_q,\Psi^1_q)_{\mR^2}\,\rho^3_{\s} 
-\bigl\{(\Psi^1_{\s},\Psi^1_q)_{\mR^2}
+|\Psi^1_q|^2\,\rho^1_{\s}\bigr\}
\frac{J^3(\bu^3)}{J^1(\bu^1)}\cos\theta^2, \\[0.2cm]
b^{33}(\bu^{13})=(\Psi^3_q,\Psi^1_{\s})_{\mR^2}
+(\Psi^3_q,\Psi^1_q)_{\mR^2}\,\rho^1_{\s} 
-\bigl\{(\Psi^3_{\s},\Psi^3_q)_{\mR^2}
+|\Psi^3_q|^2\,\rho^3_{\s}\bigr\}
\frac{J^1(\bu^1)}{J^3(\bu^3)}\cos\theta^2.
\end{eqnarray*}
Also, we define differential operators as
$$
(B_{\partial\Omega}^{ki}(\cT^{l^i}\hspace{-3pt}\rho^i_0,\mu^i_0)
\partial_\s)_{i=1,2,3}
=
\left\{
\begin{matrix}
(b_{\partial\Omega}^1(\cT^{l^1}\hspace{-3pt}\rho^1_0,\mu^1_0)
\partial_\s,0,0)
&\mbox{for}\,\ k=1, \\[0.1cm]
(0,b_{\partial\Omega}^2(\cT^{l^2}\hspace{-3pt}\rho^2_0,\mu^2_0)
\partial_\s,0)
&\mbox{for}\,\ k=2, \\[0.1cm]
(0,0,b_{\partial\Omega}^3(\cT^{l^3}\hspace{-3pt}\rho^3_0,\mu^3_0)
\partial_\s)
&\mbox{for}\,\ k=3,
\end{matrix}
\right.
$$
where $\cT^{l^i}$ ($i=1,2,3$) is the trace operator onto $\s=l^i$, i.e. 
$\cT^{l^i}\hspace{-2pt}f=f\big|_{\s=l^i}$. 
Then we set
\begin{eqnarray*}
&&\cB_0(0\,;\partial_\s)
=(B^{ki}(\cT^0\hspace{-1pt}\br_0,
\cT^0\hspace{-1pt}\br_{0,\s},\bm_0)\partial_\s)_{k,i=1,2,3}, 
\\[0.1cm]
&&\cB_0(l^i;\partial_\s)
=(B_{\partial\Omega}^{ki}(\cT^{l^i}\hspace{-3pt}\rho_0^i,\mu_0^i)
\partial_\s)_{k,i=1,2,3}, 
\end{eqnarray*}
and also set, for $(\bar{\rho}^i,\bar{\mu}^i)\in 
C^{2+\alpha,1}(\overline{\cQ^i_{0,T}})\times C^1[0,T]$ ($i=1,2,3$),
\begin{eqnarray*}
&&G^1(t)=0, \\
&&G^j(t)
=\cT^0\Bigl[\,
b^{j1}(\bu^{1j}_0)\rho^1_{0,\s}+b^{jj}(\bu^{1j}_0)\rho^j_{0,\s}
-\frac{\partial g^{1j}}{\partial\rho^1}(\bu^{1j}_0)(\bar{\rho}^1-\rho^1_0)
-\frac{\partial g^{1j}}{\partial\rho^j}(\bu^{1j}_0)(\bar{\rho}^j-\rho^j_0)
\\[0.1cm]
&&\hspace*{2.4cm}
-\frac{\partial g^{1j}}{\partial\mu^1}(\bu^{1j}_0)(\bar{\mu}^1-\mu^1_0)
-\frac{\partial g^{1j}}{\partial\mu^j}(\bu^{1j}_0)(\bar{\mu}^j-\mu^j_0)
\\[0.1cm]
&&\hspace*{2.4cm}
-\frac1{2}\int_0^1\langle D^2g^{1j}(\eta\bar{\bu}^{1j}+(1-\eta)\bu^{1j}_0)
(\bar{\bu}^{1j}-\bu_0^{1j}),\bar{\bu}^{1j}-\bu_0^{1j}\rangle\,d\eta\,\Bigr]
\quad (j=2,3), \\[0.2cm]
&&G^i_{\partial\Omega}(t)
=\cT^{l^i}\Bigl[\,
-\bigl\{b_{\partial\Omega}^i(\bar{\rho}^i,\bar{\mu}^i)
-b_{\partial\Omega}^i(\rho^i_0,\mu^i_0)\bigr\}\bar{\rho}^i_{\sigma}
+g_{\partial\Omega}^i(\bar{\rho}^i,\bar{\mu}^i)\,\Bigr]
\quad (i=1,2,3),
\end{eqnarray*}
where $Dg^{1j}$ is the Fr\'echet derivative of $g^{1j}$ and the bracket 
$\langle\cdot,\cdot\rangle$ is the respective inner product. 
Then we have the linearization of the boundary conditions:
$$
\cB_0(0\,;\partial_\s)\br^T=\bG^T(t), \quad 
\cB_0(l^i;\partial_\s)\br^T=\bG_{\partial\Omega}^T(t)
$$
for $\bG(t)=(G^1(t),G^2(t),G^3(t))$ and  
$\bG_{\partial\Omega}(t)=(G_{\partial\Omega}^1(t),
G_{\partial\Omega}^2(t),G_{\partial\Omega}^3(t))$.


{\it Step 2.} 
Let us verify that the complementary conditions hold for 
the linearized system. We refer to Lunardi \cite{Lu} for more information on
the role of the complementary conditions. 
For that purpose, we make some preparations. 
Let $\cL_0(r,\bfi\zeta)=(a^{ij}_0)_{i,j=1,2,3}$ where 
$$
a^{ii}_0=\dfrac1{\{J^i(\bu_0^i)\}^2}\zeta^2+r, \quad
a^{ij}_0=0\,\ \mbox{for}\ i\ne j.
$$
Then we have 
$$
\det\cL_0
=\prod_{i=1}^3\left[\,\dfrac1{\{J^i(\bu_0^i)\}^2}\zeta^2+r\,\right].
$$
Setting $\hat{\cL}_0=(\hat{a}_0^{ij})=(\det\cL_0)(\cL_0)^{-1}$, 
we are led to 
$$
\hat{a}^{ii}_0
=\dis\prod_{k=1,i\ne k}^3
\left[\,\dfrac1{\{J^i(\bu_0^i)\}^2}\zeta^2+r\,\right], \quad 
\hat{a}^{ij}_0=0\,\ \mbox{for}\ i\ne j.
$$
The matrix of the boundary conditions at $\s=0$ is denoted by 
$$
(B^{ki}_0(0\,;\bfi\zeta))_{i=1,2,3}
=
\left\{
\begin{matrix}
(\gamma^1,\gamma^2,\gamma^3)\hfill 
&\mbox{for}\,\ k=1, \\[0.1cm]
(\bfi b^{21}(\cT^0\bu^{12}_0)\zeta,\bfi b^{22}(\cT^0\bu^{12}_0)\zeta,0)
&\mbox{for}\,\ k=2, \\[0.1cm]
(\bfi b^{31}(\cT^0\bu^{13}_0)\zeta,0,\bfi b^{33}(\cT^0\bu^{13}_0)\zeta)
&\mbox{for}\,\ k=3, 
\end{matrix}
\right.
$$
and at $\s=l^i$ by
$$
(B^{ki}_0(l^i;\bfi\zeta))_{i=1,2,3}
=
\left\{
\begin{matrix}
(\bfi b_{\partial\Omega}^1(\cT^{l^1}\hspace*{-3pt}\rho^1_0,\mu^1_0)
\zeta,0,0)
&\mbox{for}\,\ k=1, \\[0.1cm]
(0,\bfi b_{\partial\Omega}^2(\cT^{l^2}\hspace*{-3pt}\rho^2_0,\mu^2_0)
\zeta,0)
&\mbox{for}\,\ k=2, \\[0.1cm]
(0,0,\bfi b_{\partial\Omega}^3(\cT^{l^3}\hspace*{-3pt}\rho^3_0,\mu^3_0)
\zeta)
&\mbox{for}\,\ k=3.
\end{matrix}
\right.
$$
To verify the complementary condition  
we will show that the rows of the matrix ${\cal B}_0\hat{{\cal L}}_0$ 
are linearly independent for all $r\in{\mathbb C}$, 
$\mbox{Re}\,r>0$, modulo the polynomial 
$$
P(r,\zeta) = 
\prod_{i=1}^3\{\zeta-\zeta_0^i(r)\}
$$
where $\zeta_0^i(r)
=J^i(\bu_0^i)|r|^{\frac12}e^{\bfi(\frac{\Theta}2+\frac{\pi}2)}$ 
($i=1,2,3$). 
Here we note that $\zeta^i_0(r)$ are the roots of the polynomial 
$\det{\cal L}_0(r,\bfi\zeta)$ which have a positive imaginary part. 

First let us verify the complementary condition at $\s=0$. To determine whether or not 
the complementary condition is satisfied, we have to verify that the system 
$$
\sum_{k=1}^3\omega_k{\cal B}^{ki}_0(0\,;\bfi\zeta) 
\hat{a}^{ii}_0(\zeta)\equiv0 
\quad\mbox{mod}\quad 
P(r,\zeta)=\prod_{i=1}^3(\zeta-\zeta_0^i) \quad (i=1,2,3)
$$
has the unique solution $(\omega_1,\omega_2,\omega_3)^T=\bzero$ 
or equivalently that $(\omega_1,\omega_2,\omega_3)^T=\bzero$ is the only 
vector satisfying 
$\sum_{k=1}^3\omega_k{\cal B}^{ki}_0(0\,;\bfi\zeta)\equiv0 
\quad\mbox{mod}\quad \zeta-\zeta_0^i \quad (i=1,2,3)$. 
That is, we may investigate that $(\omega_1,\omega_2,\omega_3)^T=\bzero$ is 
the only vector satisfying 
$\sum_{k=1}^3\omega_k{\cal B}^{ki}_0(0\,;\bfi\zeta_0^i)=0 \quad (i=1,2,3)$.
Thus it suffices to show that 
\begin{equation}
\det\left(\begin{array}{ccc}
\gamma^1&\gamma^2&\gamma^3 \\[0.1cm]
b^{21}(\cT^0\bu^{12}_0)&b^{22}(\cT^0\bu^{12}_0)&0 \\[0.1cm]
b^{31}(\cT^0\bu^{13}_0)&0&b^{33}(\cT^0\bu^{13}_0)
\end{array}\right)\ne0.
\label{comp_1}
\end{equation}
Indeed, in the case $(\rho^i_0,\mu^i_0)\equiv(0,0)$ ($i=1,2,3$), we have
\begin{eqnarray*}
\det\left(\begin{array}{ccc}
\gamma^1&\gamma^2&\gamma^3 \\[0.1cm]
b^{21}(\bzero)&b^{22}(\bzero)&0 \\[0.1cm]
b^{31}(\bzero)&0&b^{33}(\bzero)
\end{array}\right) &=&\det\left(\begin{array}{ccc}
\gamma^1&\gamma^2&\gamma^3 \\[0.1cm]
\sin\theta^3&-\sin\theta^3&0 \\[0.1cm]
-\sin\theta^2&0&\sin\theta^2
\end{array}\right)\nonumber\\
&=&-(\gamma^1+\gamma^2+\gamma^3)\sin\theta^2\sin\theta^3 \ne 0.
\end{eqnarray*}
Since the determinant in (\ref{comp_1}) is continuous 
with respect to $\rho^i_0$ and $\mu^i_0$ ($i=1,2,3$), we are led to 
(\ref{comp_1}) provided that $\|\rho_0^i\|_{C^1(\cI^i)}$ and $|\mu_0^i|$ 
are small enough. 

Next let us verify the complementary conditions at $\s=l^i$. 
Similarly as in the case $\s=0$, it suffices to show that 
$$
\diag(b_{\partial\Omega}^1(\cT^{l^1}\hspace*{-3pt}\rho^1_0,\mu^1_0),
b_{\partial\Omega}^2(\cT^{l^2}\hspace*{-3pt}\rho^1_0,\mu^2_0),
b_{\partial\Omega}^3(\cT^{l^3}\hspace*{-3pt}\rho^1_0,\mu^3_0))\ne0.
$$
Indeed, in the case $(\rho^i_0,\mu^i_0)\equiv(0,0)$ ($i=1,2,3$), 
we have
\begin{eqnarray*}
\det\left\{\diag(
b_{\partial\Omega}^1(0,0),
b_{\partial\Omega}^2(0,0),
b_{\partial\Omega}^3(0,0)
)\right\} 
&=&\det\left\{\diag(
|\nabla\psi(\Phi^1_*)|,
|\nabla\psi(\Phi^2_*)|,
|\nabla\psi(\Phi^3_*)|
)\right\} \\
&=&-\prod_{i=1}^3|\nabla\psi(\Phi^i_*)|\ne 0\,.
\end{eqnarray*}
Since the determinant  is also continuous 
with respect to $\rho^i_0$ and $\mu^i_0$ ($i=1,2,3$), we conclude  
$\det\bigl\{\diag(
b_{\partial\Omega}^1(\cT^{l^1}\hspace*{-3pt}\rho^1_0,\mu^1_0),
b_{\partial\Omega}^2(\cT^{l^2}\hspace*{-3pt}\rho^1_0,\mu^2_0),
b_{\partial\Omega}^3(\cT^{l^3}\hspace*{-3pt}\rho^1_0,\mu^3_0))\bigr\}
\not=0$ provided that $\|\rho_0^i\|_{C^1(\cI^i)}$ and $|\mu_0^i|$ 
are small enough. 

{\it Step 3.} 
Let us analyze the linearized system. 
Set $X=C(\cI^1)\times C(\cI^2)\times C(\cI^3)$ and 
$Y=C^2(\cI^1)\times C^2(\cI^2)\times C^2(\cI^3)$. 
Define the realization of $\cA_0$ in $X$ with homogeneous 
boundary conditions as follows
\begin{eqnarray*}
&&D(A_0)
=\{\bphi\in Y\,|\,
\bphi,\, \cA_0\bphi\in X,\,
\cB_0(0\,;\partial_\s)\bphi=\bzero,\,\cB_0(l^i;\partial_\s)\bphi=\bzero\}, 
\\
&&A_0\bphi=\cA_0\bphi.
\end{eqnarray*}
Then we have the following lemma which, in particular, characterizes the
interpolation spaces $D_{A_0}(\beta,\infty)$. For a definition of
$D_{A_0}(\beta,\infty)$ we refer to Lunardi \cite{Lu}.

\begin{lem} \label{lem:sectorial_0}\
\begin{list}{}{\topsep=0.1cm\itemsep=0cm\leftmargin=0.2cm}
\item[(i)] The linear operator 
$A_0\,:\,D(A_0)\to X$ is sectorial. 
\item[(ii)]
The characterization of the interpolation spaces
$D_{A_0}(\alpha,\infty)$ is given as
\begin{equation}
D_{A_0}(\beta,\infty)
=\left\{\begin{array}{ll}
\{\bphi\in X_\beta\mid\sum^3_{i=1}\gamma^i\varphi^i(0)=0\}
&\mbox{if}\,\ \beta\in(0,1/2), \\[0.1cm]
\{\bphi\in X_{\beta}\,|\,
\cB_0(0\,;\partial_\s)\bphi=\bzero,\,\cB_0(l^i;\partial_\s)\bphi=\bzero\}
&\mbox{if}\,\ \beta\in(1/2,1),
\end{array}\right.
\label{characterization}
\end{equation}
\end{list}
where $X_\beta := C^{2\beta}(\cI^1)\times C^{2\beta}(\cI^2)\times C^{2\beta}(\cI^3)$.
\end{lem}
{\it Proof}. 
See \cite[Section 2]{T} or adapt the argument in \cite[Section 3.1.5]{Lu} 
for the case of systems with the estimates in \cite[Theorem 12.2]{Am}. 
\qed

\medskip\noindent
Set $A_1\bphi=\cA_1\bphi$ for $\bphi\in D(A_0)$. Then we obtain 
the following lemma. 

\begin{lem} \label{lem:sectorial}
Let $A=A_0+A_1$. Then $A\,:\,D(A_0)\to X$ is a sectorial operator. 
\end{lem}
{\it Proof}. 
According to \cite[Proposition 2.4.1(ii)]{Lu}, $A$ is a sectorial operator 
in $X$ if $A_1$ is a bounded linear 
operator from $D(A_0)$ to $D_{A_0}(\alpha/2,\infty)$. Indeed, 
by means of (\ref{characterization}) and the definition of $\ba_1$, 
we have, for $\bphi\in D(A_0)$,
\begin{eqnarray*}
\|A_1\bphi\|_{D_{A_0}(\alpha/2,\infty)}
&\le&C_0\left(\sum_{i=1}^3\|\Lambda^i(\bu^i_0)\|_{C^{\alpha}(\cI^i)}\right)
\|\ba_1(\cT^0\hspace{-1pt}\br_0,\cT^0\hspace{-1pt}\br_{0,\s},\bm_0)\|
_{X_{\alpha/2}}
\|\cT^0\partial_\s^2\bphi\|_X \\[0.1cm]
&\le&\hat{C_0}\|\bphi\|_{D(A_0)}
\end{eqnarray*}
where $\alpha\in (0,1)$ and $C_0$, $\hat{C}_0$ are constants which depend on 
$\|\rho^i_0\|_{C^{1+\alpha}(\cI^i)}$ and $|\mu^i_0|$. 
This completes the proof. \qed 

\medskip\noindent
Using an estimate as in the proof of \cite[Proposition 2.4.1(ii)]{Lu}, 
we see that $A\,:\,D(A)=D(A_0)\to X$ is a sectorial operator such that 
$c_1\|\bphi\|_{D(A_0)}\le \|\bphi\|_{D(A)}\le c_2\|\bphi\|_{D(A_0)}$ 
for $\bphi\in D(A_0)$ and some constants $c_1,\,c_2>0$. 
Hence $D_A(\alpha,\infty)=D_{A_0}(\alpha,\infty)$ 
with equivalence of the respective norms. 

By virtue of  Lemma \ref{lem:sectorial}, we find that $A$ generates 
the analytic semigroup $e^{tA}$. Then we are led to the following 
proposition guaranteeing the existence of a unique solution for 
our linearized system. 

\begin{pr} \label{pr:linear}
Let us assume that $\rho^i_0\in C^{2+\alpha}(\cI^i)$ $(i=1,2,3)$ satisfy 
the compatibility conditions 
$$
\cB_0(0\,;\partial_\s)\br_0^T=\bG^T(0), \quad 
\cB_0(l^i;\partial_\s)\br_0^T=\bG_{\partial\Omega}^T(0).
$$
For $(\bar{\rho}^i,\bar{\mu}^i)\in C^{2+\alpha,1}(\overline{\cQ^i_{0,T}})
\times C^1[0,T]$ ($i=1,2,3$), the linearized system
\begin{equation}
\left\{\begin{array}{l}
\partial_t\br^T=\cA\br^T+\bF^T(\s,t), \\[0.1cm]
\cB_0(0\,;\partial_\s)\br^T=\bG^T(t), \quad
\cB_0(l^i;\partial_\s)\br^T=\bG_{\partial\Omega}^T(t), \\[0.1cm]
\rho^i_0(\cdot,0)=\rho^i_0 \quad (i=1,2,3)
\end{array}\right.
\label{L}
\end{equation}
with the notation $\cA=\cA_0+\cA_1$  has a unique solution 
such that 
\begin{eqnarray}
&&\hspace*{-0.8cm}
\sum_{i=1}^3\|\rho^i\|_{C^{2+\alpha,1}(\overline{\cQ^i_{0,T}})}
\nonumber \\[0.1cm]
&&\hspace*{-0.8cm}
\le C\sum_{i=1}^3\bigl(\|\rho^i_0\|_{C^{2+\alpha}(\cI^i)}
+\|F^i\|_{C^{\alpha,0}(\overline{\cQ^i_{0,T}})}
+\|G^i\|_{C^{(1+\alpha)/2}[0,T]}
+\|G_{\partial\Omega}^i\|_{C^{(1+\alpha)/2}[0,T]}\bigr).
\label{L-est}
\end{eqnarray}
\end{pr}
{\it Proof}. 
Adapt the argument in the proof of \cite[Theorem 5.1.19]{Lu} 
to our linearized system (\ref{L}). \qed

\medskip
Now we are ready to prove Theorem \ref{th:l-exist} by using the contraction 
principle. 

\medskip\noindent
{\it Proof of Theorem \ref{th:l-exist}.}  Set
$$
\begin{array}{rl}
{\cal D}=\bigl\{\hspace*{-0.35cm}
&(\br,\bm)\in C^{2+\alpha,1}(\overline{\cQ^1_{0,T}})\times 
C^{2+\alpha,1}(\overline{\cQ^2_{0,T}})\times
C^{2+\alpha,1}(\overline{\cQ^3_{0,T}})\times
\bigl[C^1[0,T]\,\bigr]^3\,\big|\, \\[0.2cm]
&(\rho^i(\cdot,0),\mu^i(0))=(\rho_0^i,\mu_0^i)\,\ (i=1,2,3),\,\ 
\|\br\|_{C^{2+\alpha,1}_T}+\|\bm\|_{C^1_T}\le K \bigr\}
\end{array}
$$
for some bounded positive parameters $K$ and $T$ where 
$$
\|\br\|_{C^{2+\alpha,1}_T}
=\sum_{i=1}^3\|\rho^i\|_{C^{2+\alpha,1}(\overline{\cQ^i_{0,T}})}, \quad
\|\bm\|_{C^1_T}=\sum_{i=1}^3\|\mu^i\|_{C^1[0,T]}. 
$$
Then, for $(\bar{\br},\bar{\bm})\in {\cal D}$, we define the mapping
$$
{\cal F}\,:\,{\cal D}\ni(\bar{\br},\bar{\bm})\mapsto(\br,\bm)
$$
where $\br$ is the solution of (\ref{L}) and  $\bm$ is given by 
$\bm^T=Q\,\cT^0\br^T$ for such solution $\br$. 
Once we prove that the mapping ${\cal F}$ is a contraction on ${\cal D}$ 
for suitable $K$ and $T$, the mapping ${\cal F}$ has a unique fixed point 
in ${\cal D}$ which implies that the nonlinear problem (\ref{NL}) admits 
a unique solution in $[0,T]$. 

Let us first prove that ${\cal F}$ maps ${\cal D}$ into itself. 
Note that the lower order terms in $F^i$ and $G_{\partial\Omega}^i$
can be rewritten as 
\begin{eqnarray*}
&&f^i(\bar{\bu})=f^i(\bu_0)+\int_0^1
\langle Df^i(\eta\bar{\bu}+(1-\eta)\bu_0), \bar{\bu}-\bu_0\rangle
\,d\eta, \\
&&g^i_{\partial\Omega}(\bar{\bu}^i_1)
=g^i_{\partial\Omega}(\bu^i_{1,0})+\int_0^1
\langle Dg^i_{\partial\Omega}(\eta\bar{\bu}^i_1+(1-\eta)\bu^i_{1,0}),
\bar{\bu}^i_1-\bu^i_{1,0}\rangle
\,d\eta.
\end{eqnarray*}
Here $\bar{\bu}=(\bar{\rho}^i,\bar{\rho}^i_\s,
\cT^0\hspace{-1pt}\bar{\br},\cT^0\hspace{-1pt}\bar{\br}_\s,\bar{\bm})$, 
$\bu_0=(\rho^i_0,\rho^i_{0,\s},
\cT^0\hspace{-1pt}\br_0,\cT^0\hspace{-1pt}\br_{0,\s},\bm_0)$, 
$\bar{\bu}^i_1=(\bar{\rho}^i,\bar{\mu}^i)$, 
and $\bu_{1,0}^i=(\rho_0^i,\mu_0^i)$. Moreover, $Df^i$ and 
$Dg^i_{\partial\Omega}$ are the Fr\'echet derivative of 
$f^i$ and $g^i_{\partial\Omega}$, respectively, and  
the bracket $\langle\cdot,\cdot\rangle$ is the respective inner product. 
Then, by means of (\ref{L-est}) and $\bm^T=Q\,\cT^0\br^T$, we have
\begin{eqnarray*}
\|\br\|_{C^{2+\alpha,1}_T}+\|\bm\|_{C^1_T}
&\le& \hat{C}\biggl\{\sum_{i=1}^3\bigl(\|\rho^i_0\|_{C^{2+\alpha}(\cI^i)}
+\|f^i(\bu_0)\|_{C^{\alpha}(\cI^i)}
+|g^i_{\partial\Omega}(\cT^{l^i}\bu^i_0)|\bigr) \\
&&+\sum_{j=2}^3\bigl(|b^{j1}(\cT^0\bu_0^{1j})||\cT^0\rho_{0,\s}^1|
+|b^{jj}(\cT^0\bu_0^{1j})||\cT^0\rho_{0,\s}^j|\bigr)\biggr\}+C_KT^{\nu}
\end{eqnarray*}
for $\nu=\min\{\alpha/2,\,(1-\alpha)/2\}$. Thus, choosing 
\begin{eqnarray}
K&=&2\hat{C}\biggl\{\sum_{i=1}^3\bigl(\|\rho^i_0\|_{C^{2+\alpha}(\cI^i)}
+\|f^i(\bu_0)\|_{C^{\alpha}(\cI^i)}
+|g^i_{\partial\Omega}(\cT^{l^i}\bu^i_0)|\bigr) \nonumber \\
&&+\sum_{j=2}^3\bigl(|b^{j1}(\cT^0\bu_0^{1j})||\cT^0\rho_{0,\s}^1|
+|b^{jj}(\cT^0\bu_0^{1j})||\cT^0\rho_{0,\s}^j|\bigr)\biggr\}, 
\label{K-def}
\end{eqnarray}
we conclude that there exists a time $T_1>0$ such that
\begin{equation}
\|\br\|_{C^{2+\alpha,1}_T}+\|\bm\|_{C^1_T}\le K \quad\mbox{for}\,\ T\le T_1.
\label{L-est_1}
\end{equation}
That is, ${\cal F}$ maps ${\cal D}$ into itself. 

Let us prove that the mapping ${\cal F}$ is a contraction. 
For $(\bar{\br}_1,\bar{\bm}_1),\ (\bar{\br}_2,\bar{\bm}_2)\in{\cal D}$ 
with $T\le T_1$, let 
$$
(\br_1,\bm_1)={\cal F}(\bar{\br}_1,\bar{\bm}_1), \quad 
(\br_2,\bm_2)={\cal F}(\bar{\br}_2,\bar{\bm}_2)
$$
be the solutions associated with the linearized problem (\ref{L}). 
Then, applying a similar argument to \cite[pp. 373-375]{BR} with 
$\bm^T=Q\,\cT^0\br^T$, we are led to
\[
\|\br_1-\br_2\|_{C^{2+\alpha,1}_T}+\|\bm_1-\bm_2\|_{C^1_T}
\le \hat{C}_KT^{\nu}
(\|\bar{\br}_1-\bar{\br}_2\|_{C^{2+\alpha,1}_T}
+\|\bar{\bm}_1-\bar{\bm}_2\|_{C^1_T})
\]
for $\nu=\min\{\alpha/2,\,(1-\alpha)/2\}$. 
Thus, ${\cal F}$ is a contraction for $T\le T_2$, which satisfies
$\hat{C}_KT_2^{\nu}\le1/2$. Consequently, choosing $T_*=\min\{T_1,T_2\},$
we find that ${\cal F}$ has a unique fixed point in ${\cal D}$ 
for $T\le T_*$, so that the nonlinear problem (\ref{NL}) has 
a unique solution in $[0,T]$ with (\ref{L-est_1}) for $T\le T_*$. 
Further, checking the details of the estimate for the linear system, 
we obtain for $t\in[0,T]$
$$
\sum_{i=1}^3\|\rho^i(\cdot,t)\|_{C^{1+\alpha}({\cal I}^i)}
\le m_0+C_KT^{\nu},
$$
where $m_0$ depends on $\|\rho^i_0\|_{C^{1+\alpha}(\cI^i)}$ and $|\mu^i_0|$. 
Then, there exists a time $T_3>0$ such that $m_0+C_KT^{\nu}\le 2m_0$ 
for $T\le T_3$. Thus, choosing $T_0=\min\{T_*,T_3\}$, we have 
$\sum_{i=1}^3\|\rho^i(\cdot,t)\|_{C^{1+\alpha}({\cal I}^i)}\le 2m_0$ 
for $t\in[0,T]$ with $T\le T_0$. 
It is possible to guarantee $2m_0\le\delta_0$ 
for sufficiently small $\|\rho^i_0\|_{C^{1+\alpha}(\cI^i)}$ and 
$|\mu^i_0|$, where $\delta_0$ is as in (\ref{c-inverse}). 
By Lemma~\ref{lem:junction}, $|\mu^i(t)|$ is estimated by 
$\sum_{i=1}^3\|\rho^i(\cdot,t)\|_{C^0({\cal I}^i)}$, so that 
$|\mu^i(t)|$ can be smaller than $\delta_0$ if 
$\|\rho^i_0\|_{C^{1+\alpha}(\cI^i)}$ and $|\mu^i_0|$ are small enough. 
This completes the proof of Theorem \ref{th:l-exist}. \qed

\section{Linearization}\label{lin}
In order to linearize the nonlinear system (\ref{NL}) around 
the stationary solution $\Gamma_*=\bigcup_{i=1}^3\Gamma_*^i$, 
we need to establish
the following properties of $\Psi$ at $(q,\mu^i)=(0,0)$. 

\begin{lem}\label{lem:bp}\
For the parameterizations $\Psi^i, i=1,2,3,$ in Section \ref{sec:param}, 
the following properties hold: \\[0.1cm]%
\begin{tabular}{cl}
(i)&\hspace*{-0.3cm}
$\Psi^i(\s,0,0)=\Phi_*^i(\s)$ and
$\Psi^i(\s,q,0)=\Phi_*^i(\s\mu_{\partial\Omega}^i(q)/l^i)
+qN_*^i(\s\mu_{\partial\Omega}^i(q)/l^i)$. \\[0.1cm]%
(ii)&\hspace*{-0.3cm}
$\Psi^i_{\s}(\s,0,0)=T_*^i$, $\Psi^i_q(\s,0,0)=N_*^i$, 
and $\Psi^i_{\mu}(\s,0,0)=(1-\s/l^i)T_*^i$. \\[0.1cm]%
(iii)&\hspace*{-0.3cm}
$\Psi^i_{\s\s}(\s,0,0)=(0,0)^T$, 
$\Psi^i_{\s q}(\s,0,0)=(0,0)^T$, and 
$\Psi^i_{\s\mu}(\s,0,0)=(-1/l^i)T_*^i$. \\[0.1cm]%
(iv)&\hspace*{-0.3cm}
$\Psi^i_{\s\s q}(\s,0,0)=(0,0)^T$ and 
$\Psi^i_{\s\s\mu}(\s,0,0)=(0,0)^T$. 
\end{tabular}
\end{lem}
{\it Proof.}
By the definition of $\Psi^i$, (i) is obvious. Let us prove (ii).
Differentiating $\Psi^i(\s,0,0)=\Phi_*^i(\s)$ with respect to $\s$, 
we readily derive $\Psi^i_{\s}(\s,0,0)=T_*^i(\s)$. 
Applying a similar argument to \cite{gik01}, we obtain 
$\{\mu_{\partial\Omega}^i(q)\}'|_{q=0}=0$. 
Thus (i) gives $\Psi^i_q(\s,0,0)=N_*^i(\s)$. 
Further, by the definition of $\xi^i$, we have
$$
\xi^i_{\mu}(\s,0,0)=1-\s/l^i.
$$ 
It follows from the definition of $\Psi^i$ and the Frenet-Serret formulae 
that
$$
\Psi^i_{\mu}(\s,q,\mu^i)
=\xi_{\mu}(\s,q,\mu^i)(1-q\kappa_*^i)T_*^i(\xi(\s,q,\mu^i))
=\xi_{\mu}(\s,q,\mu^i)T_*^i(\xi(\s,q,\mu^i)).
$$
Putting $(q,\mu^i)=(0,0)$, the third property of (ii) is derived. 
Finally, by using (ii), we have (iii)-(iv). \qed

\bigskip
By virtue of Lemma \ref{lem:bp}, we are led to the linearization of 
(\ref{NL}) around the stationary solution 
$\Gamma_*=\bigcup_{i=1}^3\Gamma_*^i$.

\begin{pr}\label{pr:linearization}\ 
The linearization of (\ref{NL}) is given by
\begin{equation}
\left\{\begin{array}{lcl}
\beta^i\rho^i_t=\gamma^i\rho^i_{\s\s}
&\mbox{for}&\sigma\in(0,l^i)\,, \\
\dis\sum^3_{i=1}\gamma^i\rho^i=0
&\mbox{at}&\sigma=0\,, \\
\rho^1_{\s}=\rho^2_{\s}=\rho^3_{\s}
&\mbox{at}&\sigma=0\,, \\[0.1cm]
\rho_{\s}+h^i_\ast\rho^i=0
&\mbox{at}&\sigma=l^i,
\end{array}\right.
\label{LS}
\end{equation}
for $i=1,2,3$, where $h^i_*$ is the curvature of $\partial\Omega$ at 
the point $\Phi^i_*(l^i)\in\Gamma^i_*\cap\partial\Omega$. 
\end{pr}

\noindent
We remark that (\ref{LS}) corresponds to the linearized problem which 
was derived in a formal way by Ikota and Yanagida \cite{IY}. 

\bigskip\noindent
{\it Proof of Proposition \ref{pr:linearization}}. 
Applying the same argument as in \cite[Section 3]{gik01}, using 
Lemma \ref{lem:bp} and $\kappa_*^i=0$ ($i=1,2,3$) 
we obtain from equations (\ref{rho_eq}) and the boundary conditions at $\s=l^i$, 
the first and fourth equations in (\ref{LS}). 
Thus we only derive the third equation of (\ref{LS}). To simplify 
the notation, we set 
$$
\tilde{J}^i(\rho^i,\mu^i)=J(\rho^i,\rho^i_\s,\mu^i), \quad
\tilde{g}^{1j}(\rho^1,\rho^j,\mu^1,\mu^j)
=g^{1j}(\rho^1,\rho^j,\rho^1_\s,\rho^j_\s,\mu^1,\mu^j)\,\,\ (j=2,3).
$$
Then it is easy to obtain
$$
\tilde{J}^i(\rho^i,\mu^i)=1, \quad 
\partial\tilde{J}^i(0,0)[\rho^i,\mu^i]
=\frac{d}{d\e}\tilde{J}^i(\e\rho^i,\e\mu^i)\bigg|_{\e=0}
=-\frac1{l^i}\mu^i,
$$
where $\partial\tilde{J}^i(0,0)$ is the Fr\'echet derivative of 
$\tilde{J}^i$ at $(0,0)$. Recalling the definition of 
$\tilde{g}^{1j}$ ($j=2,3$) and using  Lemma \ref{lem:bp}, we have
\begin{eqnarray*}
\partial\tilde{g}^{1j}(0,0,0,0)[\rho^1,\rho^j,\mu^1,\mu^j]
&=&-\frac1{l^1}(T_*^1,T_*^j)_{\mR^2}\mu^1
-\frac1{l^j}(T_*^1,T_*^j)_{\mR^2}\mu^j
+(T_*^1,N_*^j)_{\mR^2}\rho^j_{\s} \\[0.1cm]
&&+(N_*^1,T_*^j)_{\mR^2}\rho^1_{\s}
-\left(-\frac1{l^1}\mu^1-\frac1{l^j}\mu^j\right)\cos\theta^k
\end{eqnarray*}
for $(j,k)=(2,3)$ or $(3,2)$, where $\partial\tilde{g}^{1j}(0,0,0,0)$ 
is the Fr\'echet derivative of $\tilde{g}^{1j}$ at $(0,0,0,0)$. 
Since the angle conditions at $\s=0$ give 
\begin{eqnarray*}
&&(T_*^1,T_*^2)_{\mR^2}=\cos\theta^3, \quad 
(T_*^1,N_*^2)_{\mR^2}=-\sin\theta^3, \quad
(N_*^1,T_*^2)_{\mR^2}=\sin\theta^3, \\[0.1cm]
&&(T_*^1,T_*^3)_{\mR^2}=\cos\theta^2, \quad 
(T_*^1,N_*^3)_{\mR^2}=\sin\theta^2, \quad
(N_*^1,T_*^3)_{\mR^2}=-\sin\theta^2,
\end{eqnarray*}
it follows that
\begin{eqnarray*}
&&\partial\tilde{g}^{12}(0,0,0,0)[\rho^1,\rho^2,\mu^1,\mu^2]
=(-\rho^2_{\s}+\rho^1_{\s})\sin\theta^3, \\[0.1cm]
&&\partial\tilde{g}^{13}(0,0,0,0)[\rho^1,\rho^3,\mu^1,\mu^3]
=(\rho^3_{\s}-\rho^1_{\s})\sin\theta^2
\end{eqnarray*}
 at $\s=0$. 
Hence the linearization of $\tilde{g}^{1j}(\rho^1,\rho^j,\mu^1,\mu^j)=0$ 
($j=2,3$) is
$$
(-\rho^2_{\s}+\rho^1_{\s})\sin\theta^3=0, \quad
(\rho^3_{\s}-\rho^1_{\s})\sin\theta^2=0 \quad\mbox{at}\,\ \s=0,
$$
so that, by virtue of $\theta^j\in(0,\pi)$ ($j=2,3$), we have
$$
\rho^1_{\s}=\rho^2_{\s}=\rho^3_{\s} \quad\mbox{at}\,\ \s=0.
$$
This completes the proof. \qed

\bigskip
In \cite{IY} Ikota and Yanagida investigated linearized stability 
for the  curvature flow with a triple junction (\ref{LS}). 
They derived a criterion according to which one can determine 
whether the stationary solution is linearly stable or unstable. 
In what follows, we recall their linearized stability criterion.
The main result of \cite[Theorem 1.1]{IY} is concerned with the analysis of 
the self-adjoint eigenvalue problem associated to the linearized system 
of equations (\ref{LS}). We now recall their linearized stability result.

\begin{theo} \label{linearstability}
The maximal eigenvalue of the linearized problem (\ref{LS}),
i.e. $\beta^i=0$, $i=1,2,3$, is
negative and the stationary solution is linearly stable if and only if
one of the following conditions is fulfilled:
\begin{list}{}{\topsep=0.1cm\itemsep=0cm\leftmargin=0.5cm}
\item[(a)] either all $h^1_*$, $h^2_*$, and $h^3_*$ are positive, 
\item[(b)] or, at most one of them is non-positive, and they satisfy
\[
\gamma^1 (1+l^1 h^1_*) h^2_* h^3_* +  \gamma^2 (1+l^2 h^2_*) h^1_* h^3_* 
+ \gamma^3 (1+l^3 h^3_*) h^1_* h^2_* > 0\,.
\]
\end{list}
\end{theo}

\smallskip
We will also need a variational characterization of the linearized stability
property. To this end, let us introduce the bilinear form
\begin{equation}
I_*[\bphi, \bphi]=\sum_{i=1}^3\gamma^i\biggl\{
\int_0^{l^i}(\phi^i_{s})^2\,ds+h^i_*(\phi^i)^2|_{s=l^i}\biggr\}
\label{bilinear}
\end{equation}
for all $\bphi\in{\cal E}(\Gamma_*)$, where 
$$
{\cal E}(\Gamma_*)=\bigl\{
(\phi^1,\phi^2,\phi^3)\in H^1(0,l^1)\times H^1(0,l^2)\times H^1(0,l^3)
\,\big|\,\sum_{i=1}^3\gamma^i\phi^i(0)=0\bigr\}.
$$
This bilinear form was also considered in \cite{IY}. The following
lemma is a simple consequence of the variational characterizations of
the largest eigenvalue.  

\begin{lem} \label{linearstabilityfunctional}
Let $\lambda$ be the maximal eigenvalue of the time independent 
linearized system (\ref{LS}), i.e. we set $\beta^i=0$, $i=1,2,3$. Then 
\[
I_*[\bphi,\bphi] 
>(-\lambda)\sum_{i=1}^3\gamma^i\|\phi^i\|_{L^2(\Gamma^i_*)}^2
\quad\mbox{for}\,\,\ \bphi\in{\cal E}(\Gamma_*).
\]
\end{lem}

\begin{rem}
In order to simplify the presentation,  we will henceforth
consider only the case $\beta^i=\gamma^i$, for $i=1,2,3$. 
It is worth to note that 
the linearized stability criterion is invariant with respect 
to the positive constants $\beta^i>0$ for $i=1,2,3$.  
As it should be obvious from all the energy type estimates 
to follow,  the full nonlinear stability of the stationary solution 
will not be affected by a different choice of positive mobility 
constants $\beta^i>0$ for $i=1,2,3$.  
\end{rem}

\section{Uniqueness of the stationary solution} \label{uniqueness}
In this section we prove the uniqueness 
of the stationary solution in a small $H^2$-neighborhood. 
The inverse mapping theorem also gives a bound on 
the $H^2$-norm of $\br=(\rho^1,\rho^2,\rho^3)$ in terms of the $L^2$-norm 
of the curvature $\bkappa=(\kappa^1,\kappa^2,\kappa^3)$. 

To this end, let us introduce the function space 
\begin{equation*}
{\cal M}=\bigl\{
(\rho^1,\rho^2,\rho^3)\in H^2(0,l^1)\times H^2(0,l^2)\times H^2(0,l^3)
\,\big|\, \sum_{i=1}^3\gamma^i\rho^i(0)=0\bigr\}\,.
\end{equation*}

Then $\br$ via parametrization (\ref{param}) defines a neighboring triple junction 
configuration such that the end points lie on $\partial\Omega$.

\begin{theo}\label{theo4}
Let $I_*$ be positive. Then there exists a $H^2$-neighborhood of
$\br\equiv 0$ in ${\cal M}$, such that $\br\equiv 0$ is 
the only solution of the problem
\begin{eqnarray}
&&\kappa^i=0\,, \quad 
\sphericalangle(\partial\Omega,\Gamma^i_t)=\pi/2\,,
\label{unique-ss} \\
&&\sphericalangle(\Gamma^i(t),\Gamma^j(t))=\cos\theta^k 
\quad\mbox{for}\,\,\ i,j,k\in\{1,2,3\}\,\,\ \mbox{mutually different}\,.
\label{unique-ss2}
\end{eqnarray}
\end{theo}
{\it Proof.} The idea of the proof is to use the local inverse mapping 
theorem for the curvature operator with appropriate boundary conditions. 
The positivity of
$I_\ast$ will ensure invertibility of the linearization.

Using the notation (\ref{kappai}) of Section \ref{sec:param}  we obtain
\begin{equation*}
\kappa^i(\br) 
=a^i(\rho^i,\rho^i_\sigma,\bm(\br(0)))\rho^i_{\sigma\sigma} 
+f^i(\rho^i,\rho^i_\sigma,\bm(\br(0)))
\end{equation*}
with $a^i(\rho^i,\rho^i_\sigma,\bm(\br(0)))
=1/\{J^i(\rho^i,\rho^i_\s,\bm(\br(0)))\}^2$, 
a smooth function $f^i$, and a linear mapping $\bm$ from $\mR^3$ to $\mR^3$. 
The boundary conditions in (\ref{unique-ss}) and (\ref{unique-ss2}) can 
be written as 
$$
\begin{array}{lcl}
g^i(\rho^i(l^i),\rho^i_\sigma(l^i),\bm(\br(0)))=0
&\mbox{for}&i=1,2,3\,, \\[0.1cm]
g^i(\br(0),\br_\sigma(0),\bm(\br(0)))=0
&\mbox{for}&i=4,5
\end{array}
$$
where $g^1,\dots,g^5$ are smooth functions. We define 
${\cal N}=L^2(0,l^1)\times L^2(0,l^2)\times L^2(0,l^3)$ and 
observe that solving (\ref{unique-ss}) and (\ref{unique-ss2}) is 
equivalent to finding a zero of the mapping
$$
F\,\,:B_\e(0)\rightarrow {\cal N} \times \mR^5, \quad
F(\br) = (\kappa^1(\br),\, \kappa^2(\br),\, \kappa^3(\br),\, g^1,\dots, g^5)
$$
where $B_\e(0)$ is a ball of a radius $\e>0$ around zero 
in the space ${\cal M}$. Since $H^2$ is embedded in $C^1$ 
the mapping is well defined. Arguing similarly as in Section \ref{lin} 
we obtain
$$
\begin{array}{crcl}
DF(0)\,:&{\cal M}&\rightarrow&{\cal N}\times\mR^5 \\
&\br&\mapsto&
(\rho^1_{\sigma\sigma},\rho^2_{\sigma\sigma},\rho^3_{\sigma\sigma},
\rho^1_\sigma+h^1\rho^1,\rho^2_\sigma+h^2\rho^2, 
\rho^3_\sigma+h^3\rho^3, 
\rho^1_\sigma-\rho^2_\sigma,\rho^1_\sigma-\rho^3_\sigma)\,.
\end{array}
$$
Similarly as in \cite{IY}, 
since $I_\ast$ is positive we can conclude 
that $DF(0)$ is injective and hence the Fredholm alternative gives
that $DF(0)$ is invertible. 
Now the local inverse mapping theorem
(see e.g. \cite{Z1}) gives that there is a neighborhood around 
$\bzero$ such that only $\br\equiv \bzero$ solves (\ref{unique-ss}) 
and (\ref{unique-ss2}). \qed

\bigskip
It is worthwhile noting that the mapping $F$ analyzed in the proof 
of the above theorem is in fact a local diffeomorphism. 
Therefore its inverse mapping is locally Lipschitz continuous. 
Hence we have the following corollary. 

\begin{cor}\label{cor1}
There exist constants $C,\delta_1 >0$ such that
\[
\|\br\|_{H^2} \le C \|\bkappa\|_{L^2}
\]
provided that $\|\bkappa\|_{L^2}^2<\delta_1$ and $\br\in\mathcal{M}$
fulfills $F(\br)=(\kappa^1,\kappa^2,\kappa^3,0,0,0,0,0)$.
\end{cor}

In other words, by means of the above theorem and its corollary, 
we obtained a bound on the $H^2$-norm of the solution $\br$ 
in terms of the $L^2$-norm of the curvature 
$\bkappa=(\kappa^1,\kappa^2,\kappa^3)$ in the vicinity of 
the stationary solution $\br\equiv\bzero$ provided that $\br$ fulfills
the boundary conditions. 
This useful observation will be used several times throughout 
the rest of the paper. Although we could take the standard
$L^2$-norm in the corollary above we choose the suitable weigthed $L^2$-norm
defined in (\ref{l2norm}) in the 
corollary as this will simplify the further presentation.

\section{Governing equations for the curvature and 
other geometric quantities}\label{goveq}
In order to show the global existence and the nonlinear stability
of solutions for which the bilinear form of \cite{gik01} is positive,
we apply an energy method similar to the one used in \cite{EG1} and
\cite{IKIFB}. For such a method it is important to derive evolution 
equations for the curvature.

Let $s$ be the arc-length parameter along the phase boundary $\Gamma_t$
and let $X$ be a smooth map such that $X(\cdot,t)$ is an arc-length
parameterization of $\Gamma_t$ with
$$
\Gamma^i_t=\{X^i(s,t)\,|\,s\in[0,r^i(t)]\}
$$
where $r^i$ is smooth such that $L[\Gamma^i_t]=r^i(t)$
which is the length of $\Gamma^i_t$.
Let $N^i\,(=N^i(s,t))$ be the unit normal vector of $\Gamma^i_t$.
It can be written as
$$
N^i(s,t)=\left(\begin{array}{c} 
\cos\omega^i(s,t) \\ \sin\omega^i(s,t)
\end{array}\right).
$$
Then we have
\begin{equation}
\left\{\begin{array}{l}
N^i_s=-\kappa^i T^i, \quad T^i_s=\kappa^i N^i, \\[0.1cm]
N^i_t=-\omega^i_t T^i, \quad T^i_t=\omega^i_t N^i
\end{array}\right.
\label{lem-p1}
\end{equation}
where $T^i$ is the unit tangent vector of $\Gamma^i_t$ and
$\kappa^i$ is the curvature of $\Gamma^i_t$. In addition, we define
$$
V^i=(X^i_t,N^i)_{\mR^2}, \quad v^i=(X^i_t,T^i)_{\mR^2}
$$
and hence
\begin{equation}
X^i_t=V^iN^i+v^iT^i.
\label{velocity}
\end{equation}
Differentiating (\ref{velocity}) with respect to $s$ and
using (\ref{lem-p1}), we have
\begin{eqnarray}\label{Xts}
X^i_{ts}
&=&V^i_sN^i+V^iN^i_s+v^i_sT^i+v^iT^i_s \\
&=&(V^i_s+\kappa^i v^i)N^i+(-\kappa^i V^i+v^i_s)T^i.
\end{eqnarray}

\begin{lem} \label{lem-p2}
Let $X^i$ be a smooth arc-length parameterization as above. Then 
$$
\omega^i_t=V^i_s+\kappa^i v^i, \quad
v^i_s=\kappa^i V^i\,.
$$
\end{lem}
{\it Proof.}
Since $X^i_{ts}=X^i_{st}$ and $X^i_s=T^i$, it follows from
(\ref{lem-p1}) and (\ref{velocity}) 
that
$$
\omega^i_tN^i=(V^i_s+\kappa^i v^i)N^i+(-\kappa^i V^i+v^i_s)T^i.
$$
Thus we obtain the desired results. \qed

\bigskip\noindent
As a consequence of Lemma \ref{lem-p2} we have the following lemma.

\begin{lem} \label{lem-p3}
Let $X^i$ be a smooth arc-length parameterization of $\Gamma^i_t$ as above. Then 
the curvature $\kappa^i$ ($i=1,2,3$) satisfies the evolution equation:
$$
\kappa^i_t=V^i_{ss}+(\kappa^i)^2V^i+\kappa^i_s v^i
$$
\end{lem}
{\it Proof.}
By $\omega^i_s=\kappa^i$ and Lemma \ref{lem-p2}, we obtain
$$
\kappa^i_t=\omega^i_{st}=\omega^i_{ts}
=(V^i_s+\kappa^i v^i)_s=V^i_{ss}+\kappa^i v^i_s+\kappa^i_sv^i
=V^i_{ss}+(\kappa^i)^2 V^i+\kappa^i_sv^i .
$$
This completes the proof. \qed

\bigskip\noindent
By the assumption that $\Gamma^i_t$ meets $\partial\Omega$ at the one
end point with the angle
$\pi/2$, we have
$$
\psi(X^i(r^i(t),t))=0, \quad
(\nabla\psi(X^i),N^i)_{\mR^2}=0 \quad\mbox{at}\,\ s=r^i(t).
$$
Differentiating the identity $\psi(X^i(r^i(t),t))=0$ 
with respect to $t$ and taking into account the transversality condition
$(\nabla\psi(X^i),N^i)_{\mR^2}=0$ and the governing equation 
$X^i_t = V^i N^i +v^i T^i$ we can derive the following lemma.

\begin{lem} \label{lem-p4}
At the point where $\Gamma^i_t$ ($i=1,2,3$) meets the outer boundary
$\partial\Omega$ we have
$$
v^i+(r^i)'=0 \quad\mbox{for}\quad s=r^i(t).
$$
\end{lem}

Next we derive corresponding boundary conditions at the triple 
junction point $p(t)$. It is assumed that phase boundaries $\Gamma^1$, 
$\Gamma^2$, and $\Gamma^3$ meet at the triple junction. 
Let $i,j,k \in\{1,2,3\}$ be mutually different. 
Let $\Gamma^k$ and $\gamma^k$ be the interface and the surface energy density
between phases $i$, $j$. 
Following the arguments in Bronsard and Reitich \cite{BR} 
the angles $\theta^i$ ($i=1,2,3$) of the phases at the triple junction point 
$p(t)$ fulfill Young's law
\begin{equation}
\frac{\sin\theta^1}{\gamma^1}=\frac{\sin\theta^2}{\gamma^2}
=\frac{\sin\theta^3}{\gamma^3}
\label{Younglaw}
\end{equation}
(see \cite{Y}). Young's law can be expressed as a force balance in the
following form 
\begin{equation}
\sum_{i=1}^3 \gamma^i T^i = \sum_{i=1}^3 \gamma^i N^i =  0\,.
\label{Younglaw-alt}
\end{equation}

\bigskip
Let $p(t)\in{\mathbb R}^2$ denote a triple junction. 
At the triple junction the following boundary conditions hold:
\begin{eqnarray}
&&X^1(0,t)=X^2(0,t)=X^3(0,t)\,(=p(t)),
\label{t-pos} \\[0.1cm]
&&(X^i_s,X^j_s)_{\mR^2}=\cos\theta^k \quad\mbox{at}\,\,\ p(t)
\label{t-ang}
\end{eqnarray}
for $(i,j,k)=(1,2,3),\,(2,3,1),\,(3,1,2)$. Then we obtain the
following lemma.

\begin{lem} \label{lem:tj-bc1}
At the triple junction $p(t)$ we have the following equality:
$$
\sum_{i=1}^3\gamma^iv^i=0.
$$
\end{lem}
{\it Proof.}
Differentiating (\ref{t-pos}) with respect to $t$, we obtain
\begin{equation}
\frac{dp}{dt}=X^1_t=X^2_t=X^3_t.
\label{d-tj}
\end{equation}
For $i=1,2,3,$ it holds
$$
\bigl(\,\frac{dp}{dt},T^i\bigr)_{\mR^2}=(X^i_t,T^i)_{\mR^2}=v^i.
$$
This fact and Young's law imply that
$$
\sum_{i=1}^3\gamma^iv^i
=\sum_{i=1}^3\gamma^i\bigl(\,\frac{dp}{dt},T^i\bigr)_{\mR^2}
=\bigl(\,\frac{dp}{dt},\sum_{i=1}^3\gamma^iT^i\bigr)_{\mR^2}=0.
$$
Hence the proof is complete. \qed

\bigskip
In the next lemma we derive evolution equations and boundary
conditions for the curvature.

\begin{lem} \label{lem-p5}
A smooth solution of the curvature flow equations
\begin{equation}
V^i=\kappa^i, \quad i=1,2,3, 
\label{eq-1}
\end{equation}
with the boundary conditions
\begin{equation}
\left\{\begin{array}{lcl}
\sphericalangle(\Gamma^i(t),\Gamma^j(t))=\cos\theta^k
&\mbox{for}&i,j,k \in\{1,2,3\}\,\ \mbox{mutually different},\\[0.1cm]
\sphericalangle(\partial\Omega,\Gamma^i_t)=\pi/2
&\mbox{at}&\Gamma^i_t\cap\partial\Omega, \\[0.1cm]
\partial\Gamma^i_t\subset\partial\Omega&
\end{array}\right. \label{bc-1}
\end{equation}
fulfills when expressed in the above arc-length parameterization the
evolution equations 
$$
\kappa^i_t=\kappa^i_{ss}+(\kappa^i)^3+\kappa^i_sv^i, 
\quad i=1,2,3.
$$
Furthermore, at the triple junction $p(t)$ we have 
\begin{eqnarray}
&&\sum_{i=1}^3\gamma^i\kappa^i=0, \label{tj-bc2} \\
&&\kappa^1_s+\kappa^1v^1
=\kappa^2_s+\kappa^2v^2
=\kappa^3_s+\kappa^3v^3 \label{tj-bc3}
\end{eqnarray}
and at $\Gamma^i_t\cap\partial\Omega$ the identity 
\begin{equation}
\kappa^i_s+h^i\kappa^i=0
\label{out-bc}
\end{equation}
holds. Here $h^i$ is the curvature of $\partial\Omega$ at the point
$X^i(r^i(t),t)\in\Gamma^i_t\cap\partial\Omega$.
\end{lem}
{\it Proof.}
From (\ref{d-tj}) we deduce 
$$
\bigl(\,\frac{dp}{dt},N^i\bigr)_{\mR^2}=(X^i_t,N^i)_{\mR^2}=V^i.
$$
This fact and Young's law imply that
$$
\sum_{i=1}^3\gamma^iV^i
=\sum_{i=1}^3\gamma^i\bigl(\,\frac{dp}{dt},N^i\bigr)_{\mR^2}
=\bigl(\,\frac{dp}{dt},\sum_{i=1}^3\gamma^iN^i\bigr)_{\mR^2}=0.
$$
Since $V^i=\kappa^i$, we are led to (\ref{tj-bc2}).

Differentiating (\ref{t-ang}) with respect to $t$, we obtain
$$
(X^i_{st},X^j_s)_{\mR^2}+(X^i_s,X^j_{st})_{\mR^2}=0.
$$
It follows from (\ref{Xts}) and (\ref{bc-1}) that
$$
(V^i_s+\kappa^iv^i)\sin\theta^k
+(V^j_s+\kappa^jv^j)(-\sin\theta^k)=0.
$$
By $0<\theta^k<\pi$ we have $\sin\theta^k>0$, so that
$$
V^i_s+\kappa^iv^i=V^j_s+\kappa^jv^j.
$$
Finally, Lemma \ref{lem-p3} and the fact that 
$V^i=\kappa^i$ imply (\ref{tj-bc3}). \qed

\section{A priori estimates and global existence of a smooth solution}
\label{estimates}
The purpose of this section is to derive a priori estimates 
guaranteeing global existence of a smooth solution and its convergence 
to a steady state. 
First, we derive a priori estimates for the $L^2$-norm of the curvature. 
Next we proceed with higher order energy estimates yielding a priori 
estimates for the $H^2$-norm of the curvature. As a consequence of 
these estimates we will be able to prove exponential decay of 
the $H^2$-norm of the curvature. We remark that due to the parabolic
regularization property the solution of Theorem~\ref{th:l-exist} will become
smoother for positive time such that all derivatives in the
following computations exist. 

\subsection{First order a priori estimates}
Let us define the energy functional
\vglue-0.3truecm
\[
E[\Gamma]=\sum_{i=1}^3\gamma^iL[\Gamma^i]
\]
where $\Gamma=\bigcup_{i=1}^3\Gamma^i$ and $L[\Gamma^i]$ is
the length of $\Gamma^i$. Further, throughout Section \ref{estimates}, 
we use the following notations:
\begin{eqnarray}
&&\|\phi^i\|_{L^p}=\left(\int_{\Gamma^i_t}|\phi^i|^p\,ds\right)^{\frac1p}\ 
(1\le p<\infty), \quad 
\|\phi^i\|_{L^{\infty}}={\rm ess}\sup_{\Gamma^i_t}|\phi^i|, \\
&&\|\bphi\|_{L^p}
=\left(\sum_{i=1}^3\gamma^i\|\phi^i\|_{L^p}^p\right)^{\frac1p}, 
\quad \|\bphi\|_{L^{\infty}}=\max_i\|\phi^i\|_{L^{\infty}}
\label{l2norm}
\end{eqnarray}
for a vector function $\bphi=(\phi^1,\phi^2,\phi^3)$. 

\begin{lem} \label{12}
A smooth solution of (\ref{eq-1})-(\ref{bc-1}) fulfills the following
energy type identities 
\vspace*{-0.1cm}
\begin{flushleft}
\begin{tabular}{cl}
(i)&
$\dis\frac{d}{dt}E[\Gamma_t]+\|\bkappa\|_{L^2}^2=0,$\\[0.2cm]
(ii)&
$\dis\frac{d}{dt}\|\bkappa\|_{L^2}^2
=-2\sum_{i=1}^3\gamma^i\biggl\{\int_{\Gamma^i_t}(V^i_s)^2\,ds
+h^i(V^i)^2\big|_{s=r^i}\biggr\}
+\|\bkappa\|_{L^4}^4
+\sum_{i=1}^3\gamma^i(\kappa^i)^2v^i\big|_{s=0}$ 
\\[0.5cm]
&where $h^i$ is evaluated at $X^i(r^i(t),t)$.
\end{tabular}
\end{flushleft}
\end{lem}
{\it Proof.}
By means of the identity $L[\Gamma^i_t]=r^i(t)$, Lemmata \ref{lem-p2},  \ref{lem-p4} and \ref{lem:tj-bc1}, we have
\begin{eqnarray*}
\frac{d}{dt}E[\Gamma_t]
&=&\sum_{i=1}^3\gamma^i(r^i)'
=-\sum_{i=1}^3\gamma^i\{v^i(r^i,t)-v^i(0,t)\} \\
&=&-\sum_{i=1}^3\gamma^i\int_{\Gamma^i_t}v^i_s\,ds
=-\sum_{i=1}^3\gamma^i\int_{\Gamma^i_t}\kappa^i V^i\,ds
=-\sum_{i=1}^3\gamma^i\int_{\Gamma^i_t}(\kappa^i)^2\,ds.
\end{eqnarray*}
In order to prove (ii), we compute
\begin{equation}
2\int_{\Gamma^i_t}\kappa^i\kappa^i_t\,ds
=2\int_{\Gamma^i_t}\kappa^i\{\kappa^i_{ss}+(\kappa^i)^3+\kappa^i_sv^i\}\,ds.
\label{int01}
\end{equation}
The identity 
$$
\frac{d}{dt}\int_{\Gamma^i_t}(\kappa^i)^2\,ds
=2\int_{\Gamma^i_t}\kappa^i\kappa^i_t\,ds
+(\kappa^i)^2\big|_{s=r^i}(r^i)',
$$
and Lemma \ref{lem-p4} imply 
\begin{equation}
2\int_{\Gamma^i_t}\kappa^i\kappa^i_t\,ds
=\frac{d}{dt}\int_{\Gamma^i_t}(\kappa^i)^2\,ds
+(\kappa^i)^2v^i\big|_{s=r^i}.
\label{LHS}
\end{equation}
For the right hand side of equation (\ref{int01}) we can use the
boundary condition (\ref{out-bc}). Integration by parts yields
\begin{eqnarray}
&&\hspace*{-1cm}
2\int_{\Gamma^i_t}\kappa^i
\{\kappa^i_{ss}+(\kappa^i)^3+\kappa^i_sv^i\}\,ds \nonumber \\
&&\hspace*{-1cm}
=2\biggl\{\bigl[\kappa^i\kappa_s\bigr]_{s=0}^{s=r^i}
-\int_{\Gamma^i_t}(\kappa^i_s)^2\biggr\}\,ds
+2\int_{\Gamma^i_t}(\kappa^i)^4\,ds
+2\int_{\Gamma^i_t}\kappa^i\kappa^i_sv^i\,ds \nonumber \\
&&\hspace*{-1cm}
=-2\biggl\{\int_{\Gamma^i_t}(\kappa^i_s)^2\,ds
+h^i(\kappa^i)^2\big|_{s=r^i}\biggr\}
-2\kappa^i\kappa^i_s\big|_{s=0}+2\int_{\Gamma^i_t}(\kappa^i)^4\,ds
+2\int_{\Gamma^i_t}\kappa^i\kappa^i_sv^i\,ds.
\label{RHS-1}
\end{eqnarray}
Using the identities $v^i_s=\kappa^iV^i=(\kappa^i)^2$ and integrating by parts we obtain
\begin{eqnarray*}
\int_{\Gamma^i_t}\kappa^i\kappa^i_sv^i\,ds
&=&\bigl[(\kappa^i)^2v^i\bigr]_{s=0}^{s=r^i}
-\int_{\Gamma^i_t}\kappa^i(\kappa^i_sv^i+\kappa^iv^i_s)\,ds \\
&=&\bigl[(\kappa^i)^2v^i\big]_{s=0}^{s=r^i}
-\int_{\Gamma^i_t}\kappa^i\kappa^i_sv^i\,ds
-\int_{\Gamma^i_t}(\kappa^i)^4\,ds.
\end{eqnarray*}
Thus we have
\begin{equation}
\int_{\Gamma^i_t}\kappa^i\kappa^i_sv^i\,ds
=\frac12\biggl\{\bigl[(\kappa^i)^2v^i\bigr]_{s=0}^{s=r^i}
-\int_{\Gamma^i_t}(\kappa^i)^4\,ds\biggr\}.
\label{RHS-2}
\end{equation}
It follows from (\ref{LHS}), (\ref{RHS-1}), and (\ref{RHS-2}) that
\begin{eqnarray*}
\frac{d}{dt}\int_{\Gamma^i_t}(\kappa^i)^2\,ds
&=&-2\biggl\{\int_{\Gamma^i_t}(\kappa^i_s)^2\,ds
+h^i(\kappa^i)^2\big|_{s=r^i}\biggr\}
+\int_{\Gamma^i_t}(\kappa^i)^4\,ds \\[0.2cm]
&&-2\kappa^i\kappa^i_s\big|_{s=0}
-(\kappa^i)^2v^i\big|_{s=0}.
\end{eqnarray*}
By  (\ref{tj-bc2}) and (\ref{tj-bc3}) we have
\begin{eqnarray*}
\sum_{i=1}^3\gamma^i\Bigl\{2\kappa^i\kappa^i_s\big|_{s=0}
+(\kappa^i)^2v^i\big|_{s=0}\Bigr\}
&=&2\sum_{i=1}^3\gamma^i\kappa^i(\kappa^i_s+\kappa^iv^i)\big|_{s=0}
-\sum_{i=1}^3\gamma^i(\kappa^i)^2v^i\big|_{s=0} \\
&=&-\sum_{i=1}^3\gamma^i(\kappa^i)^2v^i\big|_{s=0}.
\end{eqnarray*}
Thus we are led to the identity
\begin{eqnarray*}
\frac{d}{dt}\sum_{i=1}^3\gamma^i\int_{\Gamma^i_t}(\kappa^i)^2\,ds
&=&-2\sum_{i=1}^3\gamma^i\biggl\{\int_{\Gamma^i_t}(\kappa^i_s)^2\,ds
+h^i(\kappa^i)^2\big|_{s=r^i}\biggr\}
+\sum_{i=1}^3\gamma^i\int_{\Gamma^i_t}(\kappa^i)^4\,ds \\
&&+\sum_{i=1}^3\gamma^i(\kappa^i)^2v^i\big|_{s=0}.
\end{eqnarray*}
Since $V^i=\kappa^i$ the proof of the lemma follows. \qed

\bigskip\noindent
Let us define a bilinear form $I$ as
$$
I[\bphi,\bphi]=\sum_{i=1}^3\gamma^i\biggl\{
\int_0^{r^i}(\phi^i_s)^2\,ds+h^i(\phi^i)^2|_{s=r^i}\biggr\}
$$
for $\bphi\in{\cal E}(\Gamma)$ where
$$
{\cal E}(\Gamma)=\bigl\{
(\phi^1,\phi^2,\phi^3)\in H^1(0,r^1)\times H^1(0,r^2)\times H^1(0,r^3)
\,\big|\,\sum_{i=1}^3\gamma^i\phi^i(0)=0\bigr\}.
$$
Since $\bV\in{\cal E}(\Gamma)$ we can rewrite the statement ii) of 
Lemma~\ref{12} as
\begin{equation}
\frac{d}{dt}\|\bkappa\|_{L^2}^2+2I[\bV,\bV]
=\|\bkappa\|_{L^4}^4
+\sum_{i=1}^3\gamma^i(\kappa^i)^2v^i\big|_{s=0}.
\label{key-1}
\end{equation}
The following lemmata are crucial in the derivation of 
a priori estimates.

\subsection{Higher order estimates for the curvature.}
We define the averaged  curvature along the curve $\Gamma^i_t$ as
$\kappa^i_{av} = \frac{1}{L[\Gamma^i_t]}\int_{\Gamma^i_t} \kappa^i ds$. 
\begin{lem}\ \label{est-2}
The following estimates for a $C^3$-curve $\Gamma^i_t$ hold true.
\vspace*{-0.3cm}
\begin{flushleft}
\begin{tabular}{cl}
(i)&
$\dis\|\kappa^i-\kappa^i_{av}\|_{L^{\infty}}
\le (L[\Gamma^i_t])^{\frac12}\|\kappa^i_s\|_{L^2}$, 
$|\kappa^i_{av}|
\le \dfrac1{(L[\Gamma^i_t])^{\frac12}}\|\kappa^i\|_{L^2}$, \\[0.2cm]
(ii)&
$\dis\int_{\Gamma^i_t}(\kappa^i)^4\,ds
\le 2 \Bigl\{L[\Gamma^i_t] \|\kappa^i_s\|_{L^2}^2
+\dfrac1{L[\Gamma^i_t]}\|\kappa^i\|_{L^2}^2\Bigr\}
\|\kappa^i\|_{L^2}^2$, \\[0.5cm]
(iii)&
There are $m_1=m_1(L[\Gamma^i_t])$ and $m_2=m_2(1/L[\Gamma^i_t])$ 
such that \\[0.2cm]
&\hspace*{2cm}
$\bigl|(\kappa^i|_{s=0})^3\bigr|
\le m_1\|\kappa^i\|_{L^2}\|\kappa^i_s\|_{L^2}^2+m_2\|\kappa^i\|_{L^2}^3$.
\end{tabular}
\end{flushleft}
\end{lem}
{\it Proof.}
The estimates in (i) are established in a standard way and we do not
present details here. By the estimates in (i), we have
$$
\|\kappa^i\|_{L^{\infty}}
\le\|\kappa^i-\kappa^i_{av}\|_{L^{\infty}}+|\kappa^i_{av}|
\le (L[\Gamma^i_t])^{\frac12}\|\kappa^i_s\|_{L^2}
+\frac1{(L[\Gamma^i_t])^{\frac12}}\|\kappa^i\|_{L^2}.
$$
It implies that
\begin{eqnarray*}
\biggl|\int_{\Gamma^i_t}(\kappa^i)^4\,ds\,\biggr|
&\le&\|\kappa^i\|_{L^{\infty}}^2\|\kappa^i\|_{L^2}^2 \\
&\le&\biggl\{(L[\Gamma^i_t])^{\frac12}\|\kappa^i_s\|_{L^2}
+\frac1{(L[\Gamma^i_t])^{\frac12}}\|\kappa^i\|_{L^2}\biggr\}^2
\|\kappa^i\|_{L^2}^2 .
\end{eqnarray*}
The statement (ii) now follows from the elementary inequality $(a+b)^2\le
2(a^2+b^2)$. 

Let $r^i_0$ be such that $\kappa^i(r^i_0,t)-\kappa^i_{av}(t)=0$. We
then obtain
$$
(\kappa^i-\kappa^i_{av})^3|_{s=0}
=-3\int_0^{r^i_0}\kappa^i_s(\kappa^i-\kappa^i_{av})^2\,ds.
$$
This implies that
\begin{eqnarray*}
\bigl|(\kappa^i-\kappa^i_{av})^3|_{s=0}\bigr|
&\le&3\int_{\Gamma^i_t}|\kappa^i_s||\kappa^i-\kappa^i_{av}|^2\,ds \\
&\le&3\|\kappa^i-\kappa^i_{av}\|_{L^{\infty}}\|\kappa^i_s\|_{L^2}
\|\kappa^i-\kappa^i_{av}\|_{L^2} \\
&\le&3(L[\Gamma^i_t])^{\frac12}\|\kappa^i_s\|_{L^2}^2
\bigl\{\|\kappa^i\|_{L^2}+(L[\Gamma^i_t])^{\frac12}|\kappa^i_{av}|\bigr\} \\
&\le&6(L[\Gamma^i_t])^{\frac12}\|\kappa^i\|_{L^2}
\|\kappa^i_s\|_{L^2}^2.
\end{eqnarray*}
Note that
$
(\kappa^i)^3=(\kappa^i-\kappa^i_{av})^3
+3\kappa^i\kappa^i_{av}(\kappa^i-\kappa^i_{av})+(\kappa^i_{av})^3.
$
Then it follows that
\begin{eqnarray*}
\|\kappa^i\kappa^i_{av}(\kappa^i-\kappa^i_{av})\|_{L^{\infty}} 
&\le& \|\kappa^i_{av}(\kappa^i-\kappa^i_{av})^2\|_{L^{\infty}}
+\|(\kappa^i_{av})^2(\kappa^i-\kappa^i_{av})\|_{L^{\infty}} \\[0.1cm]
&\le& L[\Gamma^i_t]|\kappa^i_{av}|\|\kappa^i_s\|_{L^2}^2
+(L[\Gamma^i_t])^{\frac12}|\kappa^i_{av}|^2\|\kappa^i_s\|_{L^2} \\
&\le& (L[\Gamma^i_t])^{\frac12}\|\kappa^i\|_{L^2}
\|\kappa^i_s\|_{L^2}^2
+(L[\Gamma^i_t])^{\frac12}|\kappa^i_{av}|\cdot
\frac12\bigl(|\kappa^i_{av}|^2+\|\kappa^i_s\|_{L^2}^2\bigr) \\
&\le& (L[\Gamma^i_t])^{\frac12}\|\kappa^i\|_{L^2}
\|\kappa^i_s\|_{L^2}^2
+\frac12\|\kappa^i\|_{L^2}
\biggl(\frac1{L[\Gamma^i_t]}\|\kappa^i\|_{L^2}^2
+\|\kappa^i_s\|_{L^2}^2\biggr) \\
&=& \biggl\{(L[\Gamma^i_t])^{\frac12}+\frac12\biggr\}
\|\kappa^i\|_{L^2}\|\kappa^i_s\|_{L^2}^2
+\frac1{2L[\Gamma^i_t]}\|\kappa^i\|_{L^2}^3.
\end{eqnarray*}
Thus we have
\begin{eqnarray*}
\bigl|(\kappa^i|_{s=0})^3\bigr|
&\le&\bigl|(\kappa^i-\kappa^i_{av})^3|_{s=0}\bigr|
+3\|\kappa^i\kappa^i_{av}(\kappa^i-\kappa^i_{av})\|_{L^{\infty}}
+|\kappa_{av}|^3 \\
&\le&\biggl\{9(L[\Gamma^i_t])^{\frac12}+\frac32\biggr\}
\|\kappa^i\|_{L^2}\|\kappa^i_s\|_{L^2}^2
+\biggl\{\frac3{2L[\Gamma^i_t]}
+\frac1{(L[\Gamma^i_t])^{3/2}}\biggr\}\|\kappa^i\|_{L^2}^3.
\end{eqnarray*}
Hence the proof is complete. \qed

\begin{lem}\label{lem:V-v}
For smooth solutions of (\ref{eq-1}), (\ref{bc-1}) we have 
\[
(v^1,\,v^2,\,v^3)^T = Q (V^1,\ V^2,\ V^3)^T,
\]
for any $t\in(0,T],$
where the matrix $Q$ was defined in Lemma~\ref{lem:junction}.
\end{lem}
{\it Proof.}
At the triple junction $p(t)$, we have, for all $i,j\in\{1,2,3\},$ 
$$
\frac{dp}{dt}=v^iT^i+V^iN^i=v^jT^j+V^jN^j.
$$
Taking the inner product with $T^i$ we obtain
$$
v^i=v^j(T^i,T^j)+V^j(T^i,N^j).
$$
By (\ref{t-ang}), $(T^i,T^j)=\cos\theta^k$ and $(T^i,N^j)=-\sin\theta^k$.
Thus we derive
$$
v^i-c^kv^j=-s^kV^j.
$$
If we solve this with respect to $(v^1,v^2,v^3)^T$,
we are led to the desired result. \qed

\bigskip\noindent
By Lemma \ref{lem:V-v} and $V^i=\kappa^i$, we have
\begin{equation}
\biggl|\,\sum_{i=1}^3\gamma^i(\kappa^i)^2v^i
\big|_{s=0}\,\biggr|
\le C\biggl\{\sum_{i=1}^3\gamma^i(\kappa^i)^2\big|_{s=0}\biggr\}
\biggl(\sum_{i=1}^3|\kappa^i|\big|_{s=0}\biggr)
\le \hat{C}\sum_{i=1}^3\gamma^i(\kappa^i)^3\big|_{s=0}. 
\label{tj-bc}
\end{equation}

\subsection{Structural stability of the bilinear form}
The aim of this subsection is to show that positivity of the bilinear
form $I_\ast$ is invariant with respect to small perturbations of the curve
parametrization $\rho^i$. More precisely, we will show that the
bilinear form $I$ is positive definite provided that $I_\ast$ is
positive definite and $\br=(\rho^1,\rho^2,\rho^3)$ is small in the
$C^1$-norm. Taking into account Corollary~\ref{cor1} and the
continuity of embedding $H^2\hookrightarrow C^1$ the positive definiteness
of the bilinear form is preserved if the $L^2$-
norm of the curvature $\bkappa=(\kappa^1,\kappa^2,\kappa^3)$ is small.

According to Lemma~\ref{linearstabilityfunctional}, $I_*$ is positive 
provided the maximal eigenvalue for the linearized problem is negative. 
The following lemma is a direct consequence of 
\cite[Lemma~3.1 and Prop.~3.3]{IY}.

\begin{lem}\label{15} \ 
\begin{itemize}
\item[(i)] Let $\lambda$ be the maximal eigenvalue of the linearized problem.
For $\e>0$ there exists a $\delta>0$ such that, 
for any perturbation $|h^i-h^i_*|<\delta$ and 
$|L[\Gamma^i]-L[\Gamma^i_*]|<\delta, \ i=1,2,3,$ we have
\[
I[\bphi,\bphi]>(-\lambda-\e)\|\bphi\|_{L^2}^2
\quad\mbox{for}\,\ \bphi\in{\cal E}(\Gamma).
\]
\item[(ii)] There exists a $c>0$ such that 
$$
c\,\|\bphi_s\|_{L^2}^2\le I[\bphi,\bphi]+\|\bphi\|_{L^2}^2
\quad\mbox{for}\,\ \bphi\in{\cal E}(\Gamma).
$$
\end{itemize}
\end{lem}

Using Lemma~\ref{15}, we obtain the existence of constants 
$\delta_*>0$ and $c_*>0$ such that
\begin{equation}
I[\bphi,\bphi]
>-\frac{\,\lambda\,}2\|\bphi\|_{L^2}^2+c_*\|\bphi_s\|_{L^2}^2
\quad\mbox{for}\,\ \bphi\in{\cal E}(\Gamma_t)
\label{key-2}
\end{equation}
provided that, for $i=1,2,3,$
\begin{equation}
|h^i-h^i_*|<\delta_*, \quad |L[\Gamma^i]-L[\Gamma^i_*]|<\delta_*.
\label{delta^*}
\end{equation}

\smallskip
\begin{lem} \label{length}
We have the following estimates.
\begin{itemize}
\item[(i)]
There exist constants $\delta_2$, $C>0$ such that 
$|h^i-h^i_*|\le C\|\rho^i\|_{C^0(\cI^i)}$ provided that 
$\|\rho^i\|_{C^0(\cI^i)}<\delta_2$.
\item[(ii)]
There exist constants $\delta_3$, $C>0$ such that 
$$
|L[\Gamma^i]-L[\Gamma^i_*]|\le C \|\br\|_{C^1}
\quad\mbox{and}\quad |p|\le C \|\br\|_{C^0}
$$
provided that $\|\br\|_{C^1}<\delta_3$, where $p$ is the triple junction 
of $\Gamma=\bigcup_{i=1}^3\Gamma^i$.
\end{itemize}
Here we have denoted 
$
\|\br\|_{C^{k+\alpha}}=\sum_{i=1}^3\|\rho^i\|_{C^{k+\alpha}(\cI^i)}
\quad\mbox{for}\,\,\ k\in{\mathbb N}\cup\{0\},\ \alpha\in[0,1).
$
\end{lem}
{\it Proof.}
To prove (i), we recall that $\kappa_{\partial\Omega}(X^i)$ is represented by
\begin{equation}\label{curvb}
\kappa_{\partial\Omega}(X^i)
=-\frac1{|\nabla\psi(X^i)|}
([D^2\psi(X^i)]T_{\partial\Omega}(X^i),T_{\partial\Omega}(X^i))_{\mR^2}.
\end{equation}
Since the right hand side does not depend on derivatives of $\rho^i$, 
the mean value theorem implies the second inequality of (i). 

In order to prove (ii), we have to analyze properties of the function 
$\mu^i_{\partial\Omega}$ introduced in Section~\ref{sec:param}. 
From the definition it follows that 
$\mu^i_{\partial\Omega}(0)=l^i$ and 
$p_\ast+\mu^i_{\partial\Omega}(q)T^i_\ast+qN^i_\ast\in\partial\Omega$. 
Therefore $\psi(p_\ast+\mu^i_{\partial\Omega}(q)T^i_\ast+qN^i_\ast)=0$ 
and so we can compute the derivative of $\mu^i_{\partial\Omega}(q)$ as
\[
\frac{d}{dq}\mu^i_{\partial\Omega}(q)
=-\frac{\left(
\nabla\psi(p_\ast+\mu^i_{\partial\Omega}(q)T^i_\ast+qN^i_\ast),N^i_\ast
\right)_{\mR^2}}
{\left(
\nabla\psi(p_\ast+\mu^i_{\partial\Omega}(q)T^i_\ast+qN^i_\ast),T^i_\ast
\right)_{\mR^2}} .
\]
Since  $(\nabla\psi,N^i_\ast)_{\mR^2}=0$ on $\partial\Omega$ we obtain 
$\frac{d}{d q}\mu^i_{\partial\Omega}(0) =0$. 
Now, by taking the second derivative of $\mu^i_{\partial\Omega}$ and 
taking into account the expression for the curvature $h^i_\ast$ 
at the intersection of $\partial\Omega$ and $\Gamma'_\ast$ we obtain
$\frac{d^2}{d q^2}\mu^i_{\partial\Omega}(0)=h^i_\ast$. Thus
\[
\mu^i_{\partial\Omega}(q) = l^i + h^i_\ast q^2 + o(q^2) 
\quad\hbox{as}\,\ q\to 0\,.
\]
We recall that the parameterization $\Phi^i$ of the curve $\Gamma^i$ is given by
\[
\Phi^i(\s)=p_\ast 
+\left[\mu^i+ \frac{\s}{l^i}\{\mu^i_{\partial\Omega}(\rho^i(\s))-\mu^i\}
\right]T^i_\ast+\rho^i(\s)N^i_\ast\,.
\]
Using the above property of the function $\mu^i_{\partial\Omega}$ and 
the fact $\bm^T = Q \br(0)^T$ we obtain
\[
|L[\Gamma^i] - L[\Gamma^i_\ast]| 
= \left|\int_0^{l^i} (|\Phi^i_\s (\s)|-1)\,d\s \right|
\le C \| \br\|_{C^1}\,.
\]
Similarly, as $p = \mu^i T^i_\ast +\rho^i(0) N^i_\ast$ we obtain 
$|p|\le C\|\br\|_{C^0}$. With this all statements of the lemma have
been shown. \qed

\subsection{Exponential stabilization of the solution}
\begin{lem} \label{i-Lambda}
Let $\lambda$ be the maximal eigenvalue of the linearized problem. Assume that $\lambda$ is negative.
Then there exists a $\delta_3 >0$ such that
\begin{eqnarray*}
&&\|\bkappa(t)\|_{L^2}^2
\le e^{\lambda t/2}\|\bkappa(0)\|_{L^2}^2
\quad\mbox{for}\,\,\ t\in[0,T]\,, \\[0.1cm]
&&\int_0^T\|\bkappa_s(\tau)\|_{L^2}^2\,d\tau
\le \frac1{c_*}\|\bkappa(0)\|_{L^2}^2
\end{eqnarray*}
provided that $\|\bkappa(t)\|_{L^2}^2<\delta_4$ holds on $[0,T]$ where $c_*>0$ is a constant as in (\ref{key-2}).
\end{lem}
{\it Proof.}
According to Corollary~\ref{cor1}, 
there exists a $C>0$ such that $\|\br\|_{H^2} 
\le C \|\bkappa\|_{L^2}$ for $\|\bkappa\|_{L^2}^2<\delta_1$. 
By (\ref{key-1}) we have
$$
\frac{d}{dt}\|\bkappa\|_{L^2}^2+2I[\bV,\bV]
=\|\bkappa\|_{L^4}^4
+\sum_{i=1}^3\gamma^i(\kappa^i)^2v^i\big|_{s=0}.
$$
Let us first choose $\delta_2\in(0,\delta_*)\cap(0,\delta_1)$. Then,
it follows from Lemmata \ref{est-2} and \ref{length} and the inequalities (\ref{tj-bc}), (\ref{key-2}), and (\ref{delta^*}) that there are $C>0$ such that
$$
\frac{d}{dt}\|\bkappa\|_{L^2}^2
+(-\lambda)\|\bV\|_{L^2}^2+2c_*\|\bV_{\hspace*{-2pt}s}\|_{L^2}^2
\le C\bigl(\|\bkappa\|_{L^2}+\|\bkappa\|_{L^2}^2\bigr)
(\|\bkappa\|_{L^2}^2+\|\bkappa_s\|_{L^2}^2).
$$
Since $V^i=\kappa^i$, we are led to
\begin{align}
\frac{d}{dt}\|\bkappa\|_{L^2}^2
&+\Bigl\{(-\lambda)-C\bigl(\|\bkappa\|_{L^2}+\|\bkappa\|_{L^2}^2\bigr)\Bigr\}
\|\bkappa\|_{L^2}^2 \nonumber \\
&+\Bigl\{2c_*-C\bigl(\|\bkappa\|_{L^2}+\|\bkappa\|_{L^2}^2\bigr)\Bigr\}
\|\bkappa_s\|_{L^2}^2\le0 \label{d-kappa}.
\end{align}
Then, we choose a constant $\delta_4>0$ satisfying
$$
0<\delta_4
<\min\biggl\{1,\,\frac{-\lambda}{4C},\left(\frac{-\lambda}{4C}\right)^2,
\,\frac{c_*}{2C},\,\left(\frac{c_*}{2C}\right)^2\biggr\}.
$$
If we assume 
$\|\bkappa(t)\|_{L^2}^2<\delta_4$ for $t\in[0,T]$ we then have
\begin{equation}
\frac{d}{dt}\|\bkappa(t)\|_{L^2}^2
+\frac{(-\lambda)}{2}\|\bkappa(t)\|_{L^2}^2
+c_*\|\bkappa_s(t)\|_{L^2}^2\le0.
\label{d-Lambda}
\end{equation}
Using the Gronwall inequality we obtain the desired result. \qed

\subsection{Higher order energy inequalities}\label{higher-order}
So far we have shown the exponential decay of the $L^2$-norm of the curvature
$\kappa^i$. In order to prove stabilization of the curvature in the stronger
$C^{1+\alpha}$-norm  we need to derive higher order energy
type inequalities. These estimates will enable us to conclude convergence of the curvature to
zero in the $C^{1+\alpha}$-norm. In order to derive higher order estimates 
we differentiate the curvature equation 
(see Lemma~\ref{lem-p3}) with respect to $t$ and derive an
energy estimate for $\kappa^i_t$.
To this end, let us denote
\[
w^i = \kappa^i_t \,.
\]
Then differentiating the curvature equation 
$\kappa^i_t =\kappa^i_{ss} + (\kappa^i)^3 +v^i \kappa^i_s$ 
with respect to $t$ and taking into account the commutation relation 
$\partial_t\partial_s = \partial_s\partial_t$ we obtain
\[
w^i_t =w^i_{ss} + 3 (\kappa^i)^2  w^i +v^i w^i_s + v^i_t \kappa^i_s
\]
for $i=1,2,3$. Multiplying the above equation with $w^i$ and 
integrating over $\Gamma^i_t$ yields
\begin{eqnarray}
\frac12 \frac{d}{dt}\int_{\Gamma^i_t}(w^i)^2\,ds
&=&\int_{\Gamma^i_t}w^iw^i_{ss}\,ds
+3\int_{\Gamma^i_t} (\kappa^i)^2(w^i)^2\,ds
+\int_{\Gamma^i_t}(v^iw^iw^i_s+v^i_tw^i\kappa^i_s)\,ds \nonumber \\
&=&\bigl[w^iw^i_s\bigr]^{s=r^i}_{s=0}
-\int_{\Gamma^i_t} (w^i_s)^2\,ds
+3\int_{\Gamma^i_t}(\kappa^i)^2(w^i)^2\,ds \nonumber \\
&&+\int_{\Gamma^i_t}\left(v^iw^iw^i_s+v^i_tw^i\kappa^i_s\right)\,ds\,. 
\label{A2}
\end{eqnarray}
In what follows, we analyze the boundary term 
$\bigl[w^iw^i_s\bigr]^{s=r^i}_{s=0}$ appearing in the right hand side 
of (\ref{A2}). 
First we analyze the boundary term at the triple junction position 
$s=0$. Differentiating (\ref{tj-bc2}) with respect to $t$, 
we obtain
\begin{equation}
\sum_{i=1}^3 \gamma^i w^i =0
\label{A4}
\end{equation}
at the triple junction point $p(t)$. 
It follows from (\ref{tj-bc3}) that there exists a function $G(t)$ 
such that
\[
\kappa^i_s(0,t) + \kappa^i(0,t)v^i(0,t) = G(t)
\quad \hbox{for}\,\ t\ge0\,\ \hbox{and}\,\ i=1,2,3\,.
\]
Differentiating this equation with respect to $t$, we conclude
\[
w^i_s + w^iv^i+\kappa^iv^i_t=G^\prime(t)\,.
\]
Therefore we obtain, by using (\ref{A4}), 
\begin{eqnarray}
\sum_{i=1}^3\gamma^iw^iw^i_s
&=&G^\prime(t)\sum_{i=1}^3\gamma^iw^i
-\sum_{i=1}^3\gamma^i\bigl\{(w^i)^2v^i+\kappa^iw^iv^i_t\bigr\} \nonumber \\
&=&-\sum_{i=1}^3\gamma^i\bigl\{(w^i)^2v^i+\kappa^iw^iv^i_t\bigr\}.
\end{eqnarray}
By Lemma~\ref{lem:V-v}, 
we can express the term $v^i$ as a time independent linear combination 
of curvatures $\kappa^i$ ($i=1,2,3$) evaluated at the triple junction 
$p(t)$ and so $v^i_t$ can be expressed as a time independent linear 
combination of $\kappa^i_t=w^i$ ($i=1,2,3$). 
Therefore there exists a constant $C>0$ such that
\begin{equation}
\left|
\sum_{i=1}^3 \gamma^i w^i w^i_s \big|_{s=0}
\right| \le C \|\bw\|_{L^{\infty}}^2\|\bkappa\|_{L^{\infty}}.
\label{A5}
\end{equation}

Next we proceed with the estimation of the boundary term at the point $X^i \in
\Gamma^i \cap \partial\Omega$, i.e. we consider 
$s=r^i(t)$. Notice that $r^i$ is no longer constant and its dependence 
on time $t$ has to be taken into account. We will differentiate the boundary 
condition (\ref{out-bc})
\[
\kappa^i_s(r^i(t),t) + h^i(X^i(r^i(t),t)) \kappa^i(r^i(t),t) =0
\]
with respect to $t$. 
Since $\frac{d}{dt}X^i(r^i(t),t) = X^i_t + X^i_s (r^i)^\prime
= \kappa^i N^i + \{v^i +  (r^i)^\prime\} T^i =  \kappa^i N^i$ 
(see Lemma~\ref{lem-p4}) and  
$\kappa^i_{ss} = \kappa^i_t - (\kappa^i)^3 - v^i\kappa^i_s 
=  w^i - (\kappa^i)^3- v^i\kappa^i_s$
we obtain
\begin{eqnarray}
w^i_s+h^iw^i
&=&-\left(\kappa^i_{ss} + h^i\kappa^i_s\right)(r^i)^\prime
-(\nabla h^i,\frac{d}{dt}X^i)_{\mR^2}\kappa^i \nonumber \\
&=&\{w^i-(\kappa^i)^3+(h^i- v^i)\kappa^i_s\}v^i
-(\kappa^i)^2(\nabla h^i,N^i)_{\mR^2} \nonumber \\[0.1cm]
&=&\{w^i-(\kappa^i)^3-(h^i- v^i)h^i\kappa^i\}v^i 
-(\kappa^i)^2(\nabla h^i,N^i)_{\mR^2}. \label{A6}
\end{eqnarray}
Here we have used the equations $\kappa^i_s+h^i\kappa^i=0$ and 
$v^i+(r^i)^\prime=0$ and expressed $h_i$ by the right hand side of
(\ref{curvb}).
We now denote by $\Xi^i$ the right hand side of 
(\ref{A6}). Then we get
\begin{equation}
\sum_{i=1}^3 \gamma^i w^i w^i_s \big|_{s=r^i} =
- \sum_{i=1}^3 \gamma^i h^i (w^i)^2|_{s=r^i}
+ \sum_{i=1}^3 \gamma^i w^i \Xi^i|_{s=r^i}.
\label{A7}
\end{equation}
As the outer boundary $\partial\Omega$ is assumed to be $C^3$ smooth we
obtain that the terms $h^i$, $\nabla h^i$ are uniformly bounded. 
Hence the remainder term $\sum_{i=1}^3 \gamma^i w^i \Xi^i$ 
can be estimated as
\begin{eqnarray}
\left| \sum_{i=1}^3 \gamma^i w^i \Xi^i|_{s=r^i}\right| 
&\le& C \|\bw\|_{L^{\infty}} \bigl(
\|\bkappa\|_{L^{\infty}}^2 + \|\bw\|_{L^{\infty}} \|\bv\|_{L^{\infty}} 
+ \|\bkappa\|_{L^{\infty}}^3\|\bv\|_{L^{\infty}} \nonumber \\[-0.2cm]
&&+ \|\bv\|_{L^{\infty}}^2 \|\bkappa\|_{L^{\infty}} 
+ \|\bv\|_{L^{\infty}}\|\bkappa\|_{L^{\infty}}
\bigr).
\label{A8}
\end{eqnarray}

In order to complete our estimates we have to derive $L^\infty$-
estimates on the tangential velocity $v^i$ and its time derivative 
$v^i_t$. Since $v^i_s =(\kappa^i)^2$ we have
\begin{equation}
v^i(s,t)=v^i(0,t) + \int_0^s |\kappa^i(\zeta, t)|^2\, d\zeta.
\label{A9}
\end{equation}
By Lemma~\ref{lem:V-v}, we can express $v^i(0,t)$ as a time independent 
linear combination of $\kappa^i_t$ ($i=1,2,3$) evaluated at the triple 
junction $p(t)$. Therefore there exists a constant $C>0$ such that
\begin{equation}
\|\bv\|_{L^{\infty}} 
\le C( \|\bkappa\|_{L^{\infty}} + \|\bkappa\|_{L^2}^2)
\le C \|\bkappa\|_{L^{\infty}}
\label{A10}
\end{equation}
for $\|\bkappa\|_{L^2}\le 1$. Analogously, as $v^i_t(s,t)= v^i_t(0,t) 
+ 2\int_0^s \kappa^i(\zeta, t) \kappa^i_t(\zeta,t)\, d\zeta$ and 
$v^i_t(0,t)$ is a time independent linear combination of $\kappa^i_t=w^i$ 
($i=1,2,3$) evaluated at the triple junction position $p(t)$, 
we conclude
\begin{equation}
\|\bv_t\|_{L^{\infty}} 
\le C( \|\bw\|_{L^{\infty}} + \|\bkappa\|_{L^2} \|\bw\|_{L^2} )
\le C \|\bw\|_{L^{\infty}}
\label{A11}
\end{equation}
for $\|\bkappa\|_{L^2}\le 1$.

Summarizing we have shown the following equality
\begin{eqnarray*}
&&\frac12\frac{d}{dt}\sum_{i=1}^3\gamma^i\int_{\Gamma^i_t}(w^i)^2\,ds
+I[\bw,\bw] \\
&&=\sum_{i=1}^3 \gamma^i\biggl\{ w^i \Xi^i\big|_{s=r^i}
- w^i w^i_s \big|_{s=0}
+ 3 \int_{\Gamma^i_t} (\kappa^i)^2(w^i)^2\,ds
+ \int_{\Gamma^i_t}\left( v^i w^i w^i_s + v^i_t w^i \kappa^i_s\right)\,ds
\biggr\}
\end{eqnarray*}
and, consequently, the estimate
\begin{eqnarray}
\hspace*{-0.7cm}
&&\frac12 \frac{d}{dt} \|\bw\|_{L^2}^2 + I[\bw,\bw] \nonumber \\[0.1cm]
&&\le C\bigl\{
(\|\bkappa\|_{L^{\infty}}^2 + \|\bkappa\|_{L^{\infty}}^3 
+ \|\bkappa\|_{L^{\infty}}^4 ) \|\bw\|_{L^{\infty}} 
+ \|\bkappa\|_{L^{\infty}} \|\bw\|_{L^{\infty}}^2 
+ \|\bkappa\|_{L^{\infty}}^2 \|\bw\|_{L^2}^2 \nonumber \\[0.1cm]
&&\quad
+ \|\bkappa\|_{L^{\infty}} \|\bw\|_{L^2} \|\bw_s\|_{L^2}
+ \|\bw\|_{L^{\infty}}\|\bw\|_{L^2} \|\bkappa_s\|_{L^2}
\bigr\}\,. \label{A12}
\end{eqnarray}
The application of the above inequality will be twofold. 
At first, we utilize it in order to prove a bound on 
$\|\bkappa_{ss}(t)\|_{L^2}$ uniformly for $t\in[0,T)$ where $T>0$ 
is the maximal time of existence of a $C^{2+\alpha}$ solution $\br$. 
This implies together with Theorem~\ref{th:l-exist} the possibility 
of global continuation of the $C^{2+\alpha}$ solution $\br$ up to 
the maximal time of existence $T=+\infty$ and hence 
the global existence of a $C^{2+\alpha}$ solution 
will follow.  As a second application of the above inequality we will 
prove exponential stabilization of a solution in the $H^2$-norm 
of the curvature yielding the exponential stabilization $\rho(t)$ 
in its phase-space $C^{2+\alpha}$-norm.

To accomplish this goal, we have to establish  bounds for 
$\|\bkappa\|_{L^{\infty}}$ in terms of the norms $\|\bw\|_{L^2}$ and 
$\|\bkappa\|_{L^2}$. This can be done by
taking into account the equation $w^i \equiv \kappa^i_t = \kappa^i_{ss} 
+ (\kappa^i)^3 + v^i\kappa^i_s$. From this equation we have, for $i=1,2,3,$
\begin{eqnarray}
\|\kappa^i_{ss}\|_{L^2} 
&\le& C\left(
\|w^i\|_{L^2} + \|(\kappa^i)^3\|_{L^2} + \| v^i \kappa^i_s\|_{L^2}
\right) \nonumber \\
&\le& C\left(\|w^i\|_{L^2} + \|\kappa^i\|_{L^6}^3 
+ \|\kappa^i\|_{L^{\infty}}\|\kappa^i_s\|_{L^2}\right).
\label{A13}
\end{eqnarray}
Let us denote by $\|\cdot\|_{H^k}$ the following Sobolev norm of 
the Sobolev space $H^k=W^{k,2}$
$$
\|\bphi\|_{H^k}
=\|\bphi\|_{L^2}+\|\partial^k_s\bphi\|_{L^2}.
$$
Due to the continuity of embeddings $H^2\hookrightarrow H^1$ and
$H^2\hookrightarrow L^\infty$ and using Gagliardo-Nirenberg interpolation
inequalities (cf. \cite[Lemma 5.18 and Theorem 4.17]{A}), we infer 
the existence of a constant $C_0>0$ such that
\begin{equation}
\left\{\begin{array}{l}
\|\bkappa\|_{L^{\infty}} 
\le C_0 \|\bkappa\|_{H^2}^{\frac14}\|\bkappa\|_{L^2}^{\frac34}, \quad
\|\bkappa_s\|_{L^2} 
\le C_0 \|\bkappa\|_{H^2}^{\frac12}\|\bkappa\|_{L^2}^{\frac12}, \\
\|\bw\|_{L^{\infty}}
\le C_0 \|\bw\|_{H^1}^{\frac12}\|\bw\|_{L^2}^{\frac12}\,.
\end{array}\right.
\label{gagliardo}
\end{equation}
By the Young inequality $ab\le a^p/p + b^q/q$ with $p=4/3$, $q=4$,
we have, for any $\e >0$,
\[
\|\bkappa\|_{L^{\infty}} \|\bkappa_s\|_{L^2}
\le C^2_0\|\bkappa\|_{H^2}^{\frac34} \|\bkappa\|_{L^2}^{\frac54}
\le \e \|\bkappa\|_{H^2} + C_\e\|\bkappa\|_{L^2}^5
\]
and, analogously, $\|\bkappa\|_{L^6}^3 
\le C \|\bkappa\|_{L^{\infty}}^3 
\le C_0 C\|\bkappa\|_{H^2}^{\frac34} \|\bkappa\|_{L^2}^{\frac94}
\le \e \|\bkappa\|_{H^2} + C_\e \|\bkappa\|_{L^2}^9$. 
By taking $0<\e \ll 1$ small enough we obtain from (\ref{A13})
\[
\|\bkappa_{ss}\|_{L^2} 
\le C\left(\|\bkappa\|_{L^2}+ \|\bw\|_{L^2}\right)
\quad\hbox{for}\,\ \|\bkappa\|_{L^2}\le 1\,.
\]
Consequently,
\[
\|\bkappa_s\|_{L^2} 
\le C_0 \|\bkappa\|_{H^2}^{\frac12}\|\bkappa\|_{L^2}^{\frac12}
\le \frac{C_0}{2}\left(\|\bkappa\|_{H^2} + \|\bkappa\|_{L^2} \right)
\le C(\|\bkappa\|_{L^2} + \|\bw\|_{L^2})
\]
for $\|\bkappa\|_{L^2}\le 1$. Similarly
\[
\|\bkappa\|_{L^{\infty}} \le C( \|\bkappa\|_{L^2} + \|\bw\|_{L^2})
\]
for $\|\bkappa\|_{L^2}\le 1$ where $C>0$ is a generic positive constant. 
Due to the continuity of embedding $H^1\hookrightarrow L^\infty$ 
we have $\|\bw\|_{L^{\infty}} \le C \|\bw\|_{H^1}$.

We proceed by estimating the right hand side 
of (\ref{A12}). From (\ref{gagliardo}) we have
\[
\|\bkappa\|_{L^{\infty}}^4 
\le C_0 \|\bkappa\|_{H^2} 
\le C \left(\|\bkappa\|_{L^2}+ \|\bw\|_{L^2}\right)
\le C \left( 1 + \|\bw\|_{L^2}\right)
\]
for $\|\bkappa\|_{L^2}\le 1$. Consequently, 
by using the Young inequality, we obtain
\[
\|\bkappa\|_{L^{\infty}}^2 + \|\bkappa\|_{L^{\infty}}^3 
+ \|\bkappa\|_{L^{\infty}}^4
\le C \left( 1 + \|\bw\|_{L^2} \right)
\]
for $\|\bkappa\|_{L^2}\le 1$. From the embedding 
$\|\bw\|_{L^{\infty}}\le C \|\bw\|_{H^1}$ and Young's inequality 
it follows that
\[
(\|\bkappa\|_{L^{\infty}}^2 + \|\bkappa\|_{L^{\infty}}^3 
+ \|\bkappa\|_{L^{\infty}}^4) \|\bw\|_{L^{\infty}}
\le \e \|\bw\|_{H^1}^2 + C_\e (1 + \|\bw\|_{L^2}^2)
\]
for $\|\bkappa\|_{L^2}\le 1$. Using the Gagliardo-Nirenberg inequality 
(\ref{gagliardo}) and Young's inequality, we can estimate the second summand
in (\ref{A12}) as 
\begin{eqnarray*}
\|\bkappa\|_{L^{\infty}} \|\bw\|_{L^{\infty}}^2 
\le C \|\bkappa\|_{L^{\infty}} \|\bw\|_{H^1} \|\bw\|_{L^2}
\le \e \|\bw\|_{H^1}^2 
+ C_\e \|\bkappa\|_{L^{\infty}}^2 \|\bw\|_{L^2}^2.
\end{eqnarray*}
Then, by means of $\|\bkappa\|_{L^{\infty}}^2 \le \|\bkappa\|_{L^2}^2 
+ \|\bkappa_s\|_{L^2}^2 \le 1 + \|\bkappa_s\|_{L^2}^2$ 
for $\|\bkappa\|_{L^2}\le 1$, we have
$$
\|\bkappa\|_{L^{\infty}}^2 \|\bw\|_{L^2}^2
\le (1 + \|\bkappa_s\|_{L^2}^2)\|\bw\|_{L^2}^2,
$$
and so 
$$
\|\bkappa\|_{L^{\infty}} \|\bw\|_{L^{\infty}}^2 
\le \e \|\bw\|_{H^1}^2 
+ C_\e (1 + \|\bkappa_s\|_{L^2}^2) \|\bw\|_{L^2}^2.
$$
The remaining terms in (\ref{A12}) can be easily estimated 
with help of Young's inequality as 
\begin{eqnarray*}
&&\|\bkappa\|_{L^{\infty}} \|\bw\|_{L^2}  \|\bw_s\|_{L^2} 
\le \e \|\bw\|_{H^1}^2  
+ C_\e (1 + \|\bkappa_s\|_{L^2}^2) \|\bw\|_{L^2}^2, \\[0.1cm]
&&\|\bw\|_{L^{\infty}} \|\bw\|_{L^2} \|\bkappa_s\|_{L^2} 
\le \e \|\bw\|_{H^1}^2 
+ C_\e \|\bkappa_s\|_{L^2}^2 \|\bw\|_{L^2}^2
\end{eqnarray*}
for $\|\bkappa\|_{L^2}\le 1$. 
Let us introduce $\eta(t) := 1 + \|\bkappa_s(t)\|_{L^2}^2$. 
Then, by choosing $\e>0$ sufficiently small and taking into account 
the positivity of the bilinear form $I$, we obtain
\[
I[\bw,\bw] \ge \delta \|\bw\|_{H^1}^2.
\]
Therefore the function $\|\bw(t)\|_{L^2}^2$ satisfies 
the differential inequality
\begin{equation}\label{GDI}
\frac12 \frac{d}{dt} \|\bw(t)\|_{L^2}^2 \ 
\le C_1 + C_2 \eta(t) \|\bw(t)\|_{L^2}^2\,.
\end{equation}
According to Lemma~\ref{i-Lambda}, the function $\eta$ is integrable 
on the interval $(0,T)$ and 
$$
\int_0^T\eta(t)\,dt
\le T+\int_0^{+\infty}\|\bkappa_s(t)\|_{L^2}^2\,dt <\infty
$$
provided that $T<\infty$. A Gronwall lemma type of argument applied to the differential 
inequality (\ref{GDI}) yields the existence of a $C_T$, which is
monotone increasing and bounded as long as $T$ is bounded, such that
\[
\sup_{0\le t < T} \|\bw(t)\|_{L^2}^2 \le C_T < +\infty\,.
\]
By means of Lemma~\ref{i-Lambda}, we see that $\|\bkappa\|_{L^2}$ 
is small if $\|\br(0)\|_{C^{2+\alpha}}$ is small enough. 
Furthermore, $\|\br\|_{C^{1+\alpha}}$ ($0<\alpha<1/2$) is small
provided when $\|\bkappa\|_{L^2}$ is small. In addition, using 
$\|\bkappa_{ss}\|_{L^2} \le C\left(\|\bkappa\|_{L^2} 
+ \|\bw\|_{L^2}\right)$ and the fact that the norm 
$\|\br\|_{C^{2+\alpha}}$ can be estimated by $\|\bkappa\|_{H^2}$, 
we just have shown the following conclusion.

\begin{theo} \label{theo1}
The local solution of Theorem~\ref{th:l-exist} can be extended to the
time interval $[0,\infty)$ provided that $\rho^i_0$ is small enough
in the $C^{2+\alpha}$-norm.
\end{theo}

\section{Exponential stability of stationary solutions}\label{stab}
In this section we combine all the previous results to prove exponential 
stabilization of a solution to the triple junction problem which have
initial data close to a stationary stable solution.  

\begin{theo} \label{stability}
Let the assumptions of Theorem~\ref{th:l-exist} and Theorem~\ref{theo1} 
hold and let $\Gamma_*$
be such that the bilinear form $I_*$ is positive. Then there exist constants 
$C$, $\omega$, $\delta>0$ such that
\[
\|\br(t)\|_{H^2}\le C e^{-\omega t} \|\bkappa(0)\|_{L^2}
\]
for any $t\ge 0$ and $\|\bkappa(0)\|_{L^2}<\delta$.
\end{theo}
{\it Proof.}
The proof directly follows from Lemma~\ref{i-Lambda} and 
Corollary~\ref{cor1}. \qed

\medskip
Since the $H^2$-norm of $\rho$ dominates its $C^{1+\alpha}$-norm 
and the $C^{2+\alpha}$-norm majorizes $L^2$-norm of $\kappa$ 
we can state the following consequence of the previous theorem.

\begin{cor}
Under the assumptions of Theorem~\ref{stability}
there exist constants $C,\omega, \delta>0$ such that
\[
\|\br(t)\|_{C^{1+\alpha}} 
\le C e^{-\omega t}\|\br(0)\|_{C^{2+\alpha}}
\]
for any $t\ge 0$ and $\|\br(0)\|_{C^{2+\alpha}}<\delta$.
\end{cor}

\medskip
Finally, we are able to prove exponential decay in stronger norms. 
As it was already indicated in the previous section, 
we will utilize the higher energy estimate (\ref{A12}) once more 
in order to prove exponential stabilization in the $H^2$-norm of 
the curvature $\kappa$.

Recall that, for $p\ge 2,$ we have
\begin{eqnarray}
\|\bkappa\|^p_\infty \|\bw\|_{L^{\infty}}
&\le& C^p ( \|\bkappa\|_{L^2} + \|\bw\|_{L^2})^p \|\bw\|_{L^{\infty}}
\nonumber \\
&\le& C^p ( \|\bkappa\|_{L^2} + \|\bw\|_{L^2})^{p-1} 
(\|\bkappa\|_{L^2} \|\bw\|_{L^{\infty}} 
+ \|\bw\|_{L^2} \|\bw\|_{L^{\infty}} ) \\
&\le& C ( \|\bkappa\|_{L^2} + \|\bw\|_{L^2}) 
( \|\bkappa\|_{L^2}^2 \|\bw\|_{H^1}^2 ) \nonumber \\
&\le& C \|\bkappa\|_{L^2}^2 
+ C ( \|\bkappa\|_{L^2} + \|\bw\|_{L^2}) \|\bw\|_{H^1}^2
\nonumber
\label{A14}
\end{eqnarray}
provided that $\|\bkappa\|_{L^2}\le 1$ and $\|\bw\|_{L^2}\le 1$. 
Since $\|\bw\|_{L^2}\|\bw_s\|_{L^2} \le \|\bw\|_{H^1}^2$ 
and $\|\bw\|_{L^{\infty}}\|\bw\|_{L^2} \le C \|\bw\|_{H^1}^2$ 
we conclude from (\ref{A12}), (\ref{A13}), (\ref{A14})
\[
\frac12 \frac{d}{dt} \|\bw\|_{L^2}^2 + I[\bw,\bw]
\le C \|\bkappa\|_{L^2}^2 
+ C ( \|\bkappa\|_{L^2} + \|\bw\|_{L^2}) \|\bw\|_{H^1}^2
\]
for some positive constant provided that $\|\bkappa\|_{L^2}\le 1$ 
and $\|\bw\|_{L^2}\le 1$.

Similarly as in the proof of exponential decay of $\|\bkappa\|_{L^2}$ 
we use the fact that  the full Sobolev norm $\|\bw\|_{H^1}$ can 
be estimated by the  bilinear form $I(w,w)$ as follows:
\[
\delta \|\bw\|_{H^1}^2 \le I[\bw,\bw]
\]
for some positive constant $\delta>0$. 
Taking $\|\bkappa\|_{L^2}$ and $\|\bw\|_{L^2}$ sufficiently small 
such that $C ( \|\bkappa\|_{L^2} + \|\bw\|_{L^2}) \le \delta/2$ 
we end up with the inequality
\[
\frac12 \frac{d}{dt} \|\bw\|_{L^2}^2 + \frac{\delta}{2}\|\bw\|_{H^1}^2
\le C \|\bkappa\|_{L^2}^2.
\]
Defining $y = \|\bw\|_{L^2}^2$ and using $\|\bw\|_{L^2}\le \|\bw\|_{H^1}$ 
we have
\[
\frac{dy}{dt} +\delta y \le 2C M e^{-\omega t}
\]
where $M, \omega>0$ are the modulus and rate of exponential decay of
$\|\bkappa\|_{L^2}^2$, i.e. $\|\bkappa(.,t)\|_{L^2}^2\le M e^{-\omega t}$. 
Solving the above differential inequality with respect to
$y=y(t)$ we end up with the following estimate
\[
y(t)\le y(0) e^{-\delta t} + \frac{2CM}{|\omega-\delta|}\left| e^{-\omega t} 
- e^{-\delta t}\right|\,.
\]
It means that the norm $\|\bw\|_{L^2}^2$ exponentially decays 
with the rate $\min(\delta,\omega)$. Since $\|\bkappa_{ss}\|_{L^2}
\le C(\|\bkappa\|_{L^2} + \|\bw\|_{L^2})$
and the full Sobolev norm $\|\bkappa\|_{H^1}$ dominates 
$\|\bkappa\|_{L^{\infty}}$ as well as the $\|\bkappa_s\|_{L^2}$-norm 
we obtain the following convergence result:

\begin{theo} \label{higherstability}
There exist constants $\delta$, $M$, $\omega>0$ such that the solution
of Theorem~\ref{theo1} fulfills
\[
\|\bkappa(t)\|_{H^2}^2 \le M e^{-\omega t}
\]
for all $t\ge 0$ provided that $\|\bkappa(0)\|_{H^2}^2 \le \delta$.
\end{theo}





\section*{Acknowledgements}
The research of the first two authors was supported by the Regensburger
Universit\"atsstiftung Hans Vielberth, the second author was supported
by a grant of the Sumitomo Foundation and the third author was
supported by the grant ESF-EC-0206 and a DAAD project 
within the program ``Ostpartnerschaften''.

\end{document}